\documentclass{article}
\usepackage{amsfonts}
\usepackage{amsmath}
\usepackage[all]{xy}
\newtheorem{theorem}{Theorem}[section]
\newtheorem{definition}[theorem]{Definition}

\newtheorem{corollary}[theorem]{Corollary}
\newtheorem{example}[theorem]{Example}
\newtheorem{lemma}[theorem]{Lemma}

\usepackage{mathrsfs}
\usepackage{amssymb}
\usepackage{amscd}
\begin{document}

\title{ Extensions of associative and Lie algebras  via  Gr\"obner-Shirshov bases method\footnote{Supported by the NNSF of China (11171118, 11571121).}}

\author{
Yuqun Chen  and Jianjun Qiu  \\
{\small \ School of Mathematical Sciences, South China Normal
University}\\
{\small Guangzhou 510631, P. R. China}\\
{\small yqchen@scnu.edu.cn}\\
{\small jianjunqiu@126.com}
}

\date{}
 \maketitle

\begin{abstract}\noindent   Let $\mathfrak{a},\mathfrak{b},\mathfrak{e}$ be algebras over a  field $k$. Then $\mathfrak{e}$ is an extension of $\mathfrak{a}$ by $\mathfrak{b}$ if $\mathfrak{a}$ is an ideal of $\mathfrak{e}$ and $\mathfrak{b}$ is isomorphic to the quotient algebra $\mathfrak{e}/\mathfrak{a}$.  In this paper, by using
Gr\"obner-Shirshov bases theory for associative (resp. Lie)  algebras,  we give complete characterizations of  associative (resp. Lie) algebra  extensions of
$\mathfrak{a}$  by    $\mathfrak{b}$, where $\mathfrak{b}$ is presented
by generators and relations.
\vspace{4mm}
 \vspace{4mm}
\end{abstract}


\section{Introduction}

Let $A$ and $B$ be groups. By a group extension of $A$ by $B$ we shall mean a short   exact
sequence of groups and homomorphisms
$$
1\rightarrow A\rightarrow G\rightarrow B\rightarrow 1,
$$
i.e.,  $A$ is a normal subgroup of $G$ such that $B\cong G/A$. If this is the case, we also  call
$G$ an   extension of $A$ by $B$.

The classical theory of group  extension was developed  by O. Schreier \cite{scha,schb} while the homological implications of the theory were  recognized by   S. Eilenberg and S. MacLane \cite{em47}.
Some results on characterizations of  group extensions are included in   M.  Hall \cite{hall},  S. Maclane \cite{lan}, G.P. Hochschild and  J.-P.Serre \cite{hs53}, and so on. It states in
Hall's book (\cite{hall},  p. 228) that ``it is difficult to determine the identities (in $A$) leading to
conditions for an extension, where the group $B$ is presented by generators and
relations." The Hall's problem above has been solved by Yuqun Chen in a recent paper \cite{chen08} by using Gr\"{o}bner-Shirshov bases method.

Let $\mathfrak{a},\mathfrak{b},\mathfrak{e}$ be algebras over a  field $k$ and
$$
0\rightarrow \mathfrak{a} \rightarrow  \mathfrak{e}  \rightarrow \mathfrak{b} \rightarrow 0
$$
a short exact sequence. Then $\mathfrak{e}$ is called an extension of $\mathfrak{a}$ by $\mathfrak{b}$.

Inspired by Hall's problem above, one needs to find all extensions of $\mathfrak{a}$ by $\mathfrak{b}$, where $\mathfrak{b}$ is presented
by generators and relations.

There are many results on  extensions of algebras, in particular on extensions of associative algebras, Lie algebras, Lie superalgebras, etc, see, for example, \cite{mi,mi05,ch48,fp09,fp14,ho54,hs53b,lan,mack87,mo53,Shu66}. Mostly, they deal with some special cases for extensions. There is no result to determine
conditions for an extension when the algebra $\mathfrak{b}$ is presented by generators and relations.

We discover Gr\"obner-Shirshov bases method is helpful to characterise the extensions of  $\mathfrak{a}$ by   $\mathfrak{b}$ when one deals with $\mathfrak{b}$ being presented
by generators and relations.

Gr\"{o}bner bases and Gr\"{o}bner-Shirshov bases  were invented
independently by A.I. Shirshov for ideals of free (commutative,
anti-commutative) non-associative algebras \cite{Sh62a,Shir3}, free
Lie algebras \cite{Sh,Shir3} and implicitly free associative
algebras \cite{Sh,Shir3}  (see also \cite{Be78,Bo76}), by H.
Hironaka \cite{Hi64} for ideals of the power series algebras (both
formal and convergent), and by B. Buchberger \cite{Bu70} for ideals
of the polynomial algebras. Gr\"{o}bner bases and
Gr\"{o}bner-Shirshov bases theories have been proved to be very
useful in different branches of mathematics, including commutative
algebra and combinatorial algebra. It is a powerful tool to solve
the following classical problems: normal form; word problem;
conjugacy problem; rewriting system; automaton; embedding theorem;
PBW theorem;  extension; homology;  growth function; Dehn function;
complexity; etc. See, for example, the books \cite{AL, BKu94, BuCL,
BuW, CLO, Ei}, the papers \cite{Be78, b72,Bo76,BCM, DK10, gg15,  Mikhalev92}, and the
surveys \cite{BC,BC13,BCS,BC14, BFKK00, BK03}.

The paper is organized as follows. In section 2, we review Gr\"obner-Shirshov bases  theory for associative  algebras. We give 
a complete characterization of associative algebra extensions of
$\mathfrak{a}$ by   $\mathfrak{b}$, where $\mathfrak{b}$ is presented
by generators and relations. Some applications of such a characterization are given. In section 3, we review Gr\"obner-Shirshov bases  theory for Lie algebras. We give
a complete characterization of Lie algebra extensions. Some applications of such a characterization are given.

\section{Extensions of associative algebras}

\subsection{Composition-Diamond lemma for associative algebras}

In  this subsection, we review Gr\"obner-Shirshov bases   theory for associative  algebras.

Let $X$ be a set and  $X^*$    (resp.  $S(X)$) be the free monoid (resp. semigroup) on $X$.     For any $u\in X^*$, we
denote $|u|$ the length (degree) of $u$. Let   $k\langle X_+\rangle=kS(X)\ (resp.\ k\langle X\rangle=kX^*)$ be
the free  associative  algebra on $X$ without identity (resp.  with identity) over a field $k$.

A well order $>$ on $X^*$   is called monomial if for any $u,v,w\in X^*$,
$$
u>v\Rightarrow uw>vw, \ wu>wv.
$$
A classical example of monomial order  on $X^*$ is the
deg-lex order, which first compare two words by degree (length) and
then by comparing them lexicographically, where $X$ is a well-ordered set.

Now, we assume that $X^*$ is equipped with  a monomial order $>$. For
every nonzero  polynomial  $f\in
k\langle X_+\rangle$, denote $\bar{f}$  the   leading word  of $f$ with respect to the order $>$. If the
coefficient of $\bar{f}$ in $f$ is equal to 1, then  $f$ is called
monic.  If $S\subseteq k\langle X_+\rangle$ with each $s\in S$  monic, then  $S$ is also called  monic.

Let $f$ and $g$ be two monic polynomials in $k\langle X_+\rangle$.
There are two kinds of compositions:

$(i)$ If   $w=\bar{f}a=b\bar{g}$ for some
$a,b\in X^*$ with $|\bar{f}|+|\bar{g}|>|w|$, then the polynomial
 $(f,g)_w=fa-bg$ is called the intersection composition of $f$ and
$g$ with respect to the ambiguity  $w$.

$(ii)$ If  $w=\bar{f}=a\bar{g}b$ for some $a,b\in X^*$, then the
polynomial $(f,g)_w=f - agb$ is called the inclusion composition of
$f$ and $g$ with respect to the ambiguity  $w$.

Let $S\subseteq$ $k\langle X_+\rangle$ be monic. The composition $(f,g)_w$ is called trivial modulo $(S, w)$ if
$(f,g)_w=\sum\alpha_i a_i s_i b_i$, where each $\alpha_i\in k$,
$a_i,b_i\in X^{*}, \ s_i\in S$ and $a_i \overline{ s_i }b_i<w$. If
this is the case, then we write
$
(f,g)_w\equiv0\  mod(S,w).
$

A monic  set $S$  is called  a
Gr\"{o}bner-Shirshov basis  in   $k\langle X_+\rangle$ if any
composition  $(f,g)_w$ of polynomials $f, g \in S$ is trivial modulo $(S,w)$.

A Gr\"{o}bner-Shirshov basis $S$ in   $k\langle X_+\rangle$ is minimal if there is no inclusion composition in $S$. Let $\mathfrak{c}=k\langle X_+| R\rangle$ be an associative algebra
generated by  $X$  and defining relations $R$.
Then for any monomial order on $X^*$, one can find a minimal
Gr\"{o}bner-Shirshov basis $S$ in   $k\langle X_+\rangle$ such that $\mathfrak{c}=k\langle X_+| R\rangle=k\langle X_+| S\rangle$, see, \cite{BC14}.

The following theorem was proved by A.I. Shirshov \cite{Sh} for   Lie
algebras (with deg-lex order) in 1962 (see also L.A. Bokut \cite{b72}).
In 1976, L.A. Bokut \cite{Bo76} specialized the approach of A.I. Shirshov to
associative algebras (see also G.M.  Bergman \cite{Be78}). For commutative
polynomials, this lemma is known as the  Buchberger's Theorem (see
\cite{bu65,Bu70}).

\begin{theorem}\label{th2.1}   (Composition-Diamond lemma for associative algebras) \ Let $S\subseteq$ $k\langle X_+\rangle$ be monic, $Id(S)$    the ideal of $k\langle X_+\rangle$ generated by $S$ and  $>$  a monomial order on $X^*$. Then the following statements are equivalent:
\begin{enumerate}
\item[(a)] $S $ is a Gr\"{o}bner-Shirshov basis  in $k \langle X_+\rangle$.
\item[(b)] $f\in Id(S)\Rightarrow \bar{f}=a\bar{s}b$
for some $s\in S$ and $a,b\in  X^*$.
\item[(c)] $Irr(S) = \{ u \in S(X) |  u \neq a\bar{s}b ,s\in S,a ,b \in X^*\}$
is a  linear basis of the associative algebra $k\langle X_+ |
S \rangle=k\langle X_+
  \rangle/Id(S)$.
\end{enumerate}
\end{theorem}

\subsection{A characterization of extensions of associative algebras}

 Let $\mathfrak{a},\mathfrak{b}$ be associative   algebras over a  field $k$. In this subsection, by using  Gr\"{o}bner-Shirshov bases theory for associative algebras, we give a complete characterization of   extensions  of   $\mathfrak{a}$ by   $\mathfrak{b}$, where  $\mathfrak{b}$ is presented by generators and relations.

\begin{definition} Let $ \mathfrak{c}$ and $ \mathfrak{d}$ be two associative algebras. Then
$ \mathfrak{c}$  is   called a  compatible $ \mathfrak{d}$-bimodule  if
$ \mathfrak{c}$  is a $\mathfrak{d}$-bimodule and for any $c, c'\in \mathfrak{c}$ and $ d \in \mathfrak{d}$,
$$
(c\cdot d)c'=c(d\cdot c'), \ \ (d\cdot c)c'=d\cdot (cc'), \ \ (cc')\cdot d=c(c'\cdot d),
$$
 where $\cdot$ is the $\mathfrak{d}$-bimodule operation.
\end{definition}

Let $\mathfrak{a}$  be an associative algebra with  a linear basis
$A$ and the    multiplication table: $aa'=\{aa'\}, a, a'\in A$, where
$\{aa'\}$ is a linear combination of elements in $A$.  This means
that $\mathfrak{a}=k\langle A_+| aa'-\{aa'\}, a, a'\in A\rangle$.
Let $\mathfrak{b}=k\langle B_+| R\rangle$ be an associative algebra
generated by a set $B$  with defining relations $R$.

Assume that $>_{_{B^*}}$ is a monomial order on $B^*$, $>_{_{A}}$ a well order on $A$ and $>_{_{A^*}}$ the deg-lex order on  $A^*$. Let
$>_{_{(A\cup B)^*}}$ be the lex order on $(A\cup B)^*$ with respect to  the order  $>_{_{ A\cup B}}$  on $A\cup B $ such that $b>_{_{ A\cup B}}a$
  for any $a\in A, b\in B$.

Now we define a monomial order  $(A\cup B)^*$.

For any $w\in (A\cup B)^*$, $w$ can be uniquely expressed as
$$
w=u_0v_1u_1v_2u_2\cdots v_nu_n,\ \  n\geq 0,
$$
where $v_i\in S(B), i=1,2,\cdots, n$,  $u_0, u_n\in A^*$ and $u_j\in S(A),  j=1,2, \cdots, n-1$. Denote $w_{_{B}}=v_1 v_2 \cdots v_n$ for $n>0$  and
    $w_{_{B}}=1 $ for $n=0$,  $w_{_{A}}= u_0u_1u_2\cdots u_n$, and
$
wt(w)=(w_{_{B}},w_{_{A}}, w).
$

For any $w, w'\in (A\cup B)^*$, define
\begin{eqnarray}\label{or1}
w>w'\ \mbox{if }\ wt(w)>wt(w')\ \ \ \mbox{lexicographically}
\end{eqnarray}
It is clear that such an order is a monomial order  on $(A\cup B)^*$, which is an extension of the order  $>_{B^*}$ on $B^*$.

For convenience, we may assume that $R$  is a  minimal
Gr\"{o}bner-Shirshov basis in $k\langle B_+\rangle$ with respect to  $>_{B^*}$.

Let   $\lfloor\  \rfloor : R \rightarrow \mathfrak{a}$ be a map, which is called a factor set of $\mathfrak{b}$ in $\mathfrak{a}$. Suppose that $\mathfrak{a}$ is a  compatible $k\langle B_+\rangle$-bimodule with the module operation $\cdot$.  Define
$$
E_{ (\mathfrak{a}, \mathfrak{b},\lfloor\  \rfloor,\cdot)}=k\langle (A\cup B)_+|R_{(\lfloor\  \rfloor, \cdot)}\rangle,
$$
where  $R_{(\lfloor\  \rfloor, \cdot)}$  is the set consisting of the
following polynomials in $k\langle (A\cup B)_+\rangle$:
\begin{eqnarray*}
f_{aa'}:&=& aa'- \{aa'\},\ \  a,a'\in A,\\
 r_{ab}:&=&ab-a\cdot b,\ \ \ \ \ a\in A, b\in B, \\
l_{ba}:&=&ba-b\cdot a,\ \ \ \ \ a\in A, b\in B,\\
 h_f:&=&f-\lfloor f\rfloor,\ \ \ \ \ \ \ f\in R.
\end{eqnarray*}

\begin{lemma} \label{le2.3}  With the order defined by (\ref{or1}),
$R_{(\lfloor\  \rfloor, \cdot)}$ is a Gr\"{o}bner-Shirshov basis in
$k\langle (A\cup B)_+\rangle$ if and only if the following two  conditions hold:
\begin{enumerate}
\item[(a)] For any $a\in A, f\in R$,
 $
f\cdot a =\{\lfloor f\rfloor a\},\ \ a\cdot f=\{ a\lfloor f\rfloor\}
$ in $\mathfrak{a}$.
\item[(b)] For any intersection  composition $(f_1, f_2)_w$, $f_1, f_2\in R,$
$w=\overline{f_1}c=d\overline{f_2}, c,d\in B^*$  with  $|\overline{f_1}|+|\overline{f_2}|>|w|$, if
$
(f_1,f_2)_w=f_1c-df_2=\sum \alpha_i c_ig_id_i,
 $
where each $\alpha_i\in k,  c_i, d_i\in B^*, g_i\in R, c_i\overline{g_i}d_i<w$, then in $\mathfrak{a}$,
$$
\lfloor f_1\rfloor \cdot c-d\cdot\lfloor f_2\rfloor=\sum \alpha_i
c_i\cdot\lfloor   g_i \rfloor  \cdot d_i.
$$
\end{enumerate}
Moreover, if this is the case,
$
E_{ (\mathfrak{a}, \mathfrak{b},\lfloor\  \rfloor,\cdot)}
$
is an extension of $\mathfrak{a}$ by $\mathfrak{b}$.
\end{lemma}

\noindent{\bf Proof.} Since $R$ is a minimal  Gr\"{o}bner-Shirshov basis in $k\langle
B_+\rangle$, there is no inclusion composition in $R_{(\lfloor\ \rfloor,\cdot)}$. Then, the  only possible intersection compositions in
$R_{(\lfloor\ \rfloor,\cdot)}$ are related to the following
ambiguities:
$$
aa'a'',\ \overline{f}a=vba,  \ a\overline{f}=abv, \ \overline{f_1}c=d\overline{f_2}, \ aba',\ bab',\ baa',\ aa'b,
$$
where $a,a', a''\in A,\ b,b'\in B$, $f,f_1,f_2 \in R$, $v,c, d\in
B^*$, $|\overline{f_1}|+|\overline{f_2}|>|\overline{f_1}c|$.

If $w=aa'a''$, then $(f_{aa'}, f_{a'a''} )_w\equiv 0\  mod(R_{(\lfloor\  \rfloor, \cdot)},w)$.

If $w=\overline{f}a=vba$, then
\begin{eqnarray*}
 (h_f, l_{ba})_{w}&=& (f-\lfloor f\rfloor)a-v(ba-b\cdot a)
 \equiv   (f-\overline{f})a-\lfloor f\rfloor a+v(b\cdot a)\\
&\equiv&  (f-\overline{f})\cdot a-\lfloor f\rfloor a+(vb)\cdot a
 \equiv  f\cdot a -\{\lfloor f\rfloor a\}\ \ mod (R_{(\lfloor\
\rfloor,\cdot)},w).
\end{eqnarray*}
Thus, $ (h_f, l_{ba})_{w}\equiv 0\ \ mod (R_{(\lfloor\
\rfloor,\cdot)},w)$ if and only if $f\cdot a =\{\lfloor f\rfloor a\}$. \\

 Similarly,
  $(r_{ab},h_f)_{w}\equiv 0 \  mod (R_{(\lfloor\  \rfloor, \cdot)},w)$ if and only if  $a\cdot f=\{ a\lfloor f\rfloor\}$, where  $w=abv =a\overline{f}$.

If $w=\overline{f_1}c=d\overline{f_2}$  with  $|\overline{f_1}|+|\overline{f_2}|>|w|$, where $f_1,f_2\in R$, then
\begin{eqnarray*}
(h_{f_1},h_{f_2})_w=(f_1-\lfloor f_1\rfloor)c-d(f_2-\lfloor
f_2\rfloor) = (f_1c-df_2)-\lfloor f_1\rfloor c+d\lfloor f_2\rfloor.
\end{eqnarray*}
Since  $R$ is a  Gr\"{o}bner-Shirshov basis in $k\langle B_+\rangle$, we have
$$
(f_1,f_2)_w=f_1c-df_2=\sum \alpha_i c_ig_id_i,
 $$
where each $\alpha_i\in k,\ g_i\in R,\ c_i, d_i\in B^*$ and  $c_i\overline{g_i}d_i<w$.
Then,
\begin{eqnarray*}
(h_{f_1},h_{f_2})_w
&=& \sum \alpha_i c_i(g_i-\lfloor g_i \rfloor)d_i +\sum \alpha_i c_i\lfloor g_i \rfloor d_i -\lfloor f_1\rfloor c+d\lfloor f_2\rfloor\\
&\equiv& \sum \alpha_i c_i\lfloor g_i \rfloor d_i -\lfloor f_1\rfloor c+d\lfloor f_2\rfloor\\
&\equiv& \sum \alpha_i c_i\cdot\lfloor   g_i \rfloor \cdot d_i -\lfloor f_1\rfloor \cdot c+d\cdot\lfloor f_2\rfloor \ mod (R_{(\lfloor\  \rfloor, \cdot)},w).
\end{eqnarray*}
Thus, $ (h_{f_1},h_{f_2})_w\equiv 0 \  mod (R_{(\lfloor\  \rfloor, \cdot)},w)$ if and only if
$$
\lfloor f_1\rfloor \cdot c-d\cdot\lfloor f_2\rfloor=\sum \alpha_i
c_i\cdot\lfloor   g_i \rfloor \cdot d_i.
$$

Since $\mathfrak{a}$  is   a  compatible $k\langle B_+\rangle$-bimodule,  it is easy to check that
\begin{eqnarray*}
 (r_{ab},l_{ba'})_{w}&=&0 \ \ mod (R_{(\lfloor\  \rfloor,\cdot)},w),\ \ w=aba',\\
  (l_{ba},r_{ab'})_{w}&=&0 \ \ mod (R_{(\lfloor\  \rfloor,\cdot)},w),\ \ w=bab',\\
   (l_{ba},f_{aa'})_{w}&=& 0\ \ mod (R_{(\lfloor\  \rfloor,\cdot)},w),\ \ w=baa',\\
   (f_{aa'},r_{a'b})_{w}&=& 0\ \  mod (R_{(\lfloor\  \rfloor,\cdot)},w),\ \ w=aa'b.
\end{eqnarray*}
For example,
if $w=aba'$, then
\begin{eqnarray*}
 (r_{ab},l_{ba'})_{w}&=&(ab-a\cdot b)a'-a(ba'-b\cdot a')
 \equiv-(a\cdot b)a'+a(b\cdot a')\\
&\equiv& \{(a\cdot b)a'\}-\{a(b\cdot a')\}
 \equiv  0 \ \ mod (R_{(\lfloor\  \rfloor,\cdot)},w).
\end{eqnarray*}

Therefore, we can see that $R_{(\lfloor\  \rfloor, \cdot)}$ is a
Gr\"{o}bner-Shirshov basis in $k\langle (A\cup B)_+\rangle$ if and
only if the conditions  $(a)$ and $(b)$ hold.

If $R_{(\lfloor\  \rfloor, \cdot)}$ is a Gr\"{o}bner-Shirshov basis,
then by Composition-Diamond lemma for associative algebras, we
have that $Irr(R_{(\lfloor\  \rfloor, \cdot)})=A\cup Irr(R)$ is a linear
basis of $E_{(\mathfrak{a}, \mathfrak{b},\lfloor\  \rfloor,\cdot)}$,
where $ Irr(R)=\{w\in S(B)|w\neq a\bar{f}b, f\in R, a, b\in B^*\}. $
 Thus, $E_{(\mathfrak{a}, \mathfrak{b},\lfloor\  \rfloor,\cdot)}=
 \mathfrak{a}\oplus \mathfrak{b}$ as linear space with the multiplication: for any $a,a'\in A$ and $u,u'\in Irr(R)$,
$$
(a+u)(a'+u')=\{aa'\}+a \cdot u'+u\cdot a'+\sum \delta_iu_i+\sum \beta_j c_j\cdot\lfloor g_j\rfloor\cdot d_j,
$$
where $uu'=\sum \delta_iu_i+\sum \beta_j c_jg_jd_j$ in
$k\langle B_+\rangle$ and $\delta_i, \beta_j\in k$, $u_i\in Irr(R)$,
$g_j\in R, c_j, d_j\in B^*$.

It is easy to see that
$\mathfrak{a}$ is an ideal of $E_{(\mathfrak{a}, \mathfrak{b},\lfloor\
\rfloor,\cdot)}$. Define
$$
\pi_{  \llcorner\lrcorner,\cdot} : E_{(\mathfrak{a}, \mathfrak{b},\lfloor\
\rfloor,\cdot)}= \mathfrak{a}\oplus \mathfrak{b}\rightarrow \mathfrak{b},\ \tilde{a}+\tilde{b}\mapsto \tilde{b},\ \tilde{a}\in \mathfrak{a}, \tilde{b}\in \mathfrak{b}.
$$
It follows that
$
0\rightarrow \mathfrak{a} \overset{i_{ \llcorner\lrcorner,\cdot}}{\longrightarrow}  E_{(\mathfrak{a}, \mathfrak{b},\lfloor\
\rfloor,\cdot)} \overset{\pi_{ \llcorner\lrcorner,\cdot}}{\longrightarrow} \mathfrak{b}\rightarrow 0
$
is a short exact   sequence, where $ i_{ \llcorner\lrcorner,\cdot} $ is the inclusion map. Therefore, $E_{ (\mathfrak{a}, \mathfrak{b},\lfloor\  \rfloor,\cdot)}$
is an extension of $\mathfrak{a}$ by $\mathfrak{b}$.
\hfill $ \square$\\

Two extensions
  $0\rightarrow \mathfrak{a}  \overset{i}\rightarrow \mathfrak{e}\overset{\pi}\rightarrow \mathfrak{b} \rightarrow 0$ and $0\rightarrow \mathfrak{a}  \overset{i'}\rightarrow \mathfrak{e'}\overset{\pi'}\rightarrow \mathfrak{b} \rightarrow 0$
 are said to be equivalent   if there exists an algebra
homomorphism $\varphi:\mathfrak{e}\rightarrow \mathfrak{e}'$ such
that the following diagram is commutative:
$$
\xymatrix{
 0 \ar[r]&\mathfrak{a} \ar[r]^{i} \ar@{=}[d]_{}& \mathfrak{e}\ar[r]^{\pi}\ar[d]^{\varphi} &    \mathfrak{b} \ar@{=}[d]^{} \ar[r]&0  \\
       0 \ar[r] &   \mathfrak{a} \ar[r]^{i'}  & \mathfrak{e}'\ar[r]^{\pi'}& \mathfrak{b}\ar[r]  &    0       }
$$

Let  $0\rightarrow \mathfrak{a}  \overset{i}\rightarrow \mathfrak{e}\overset{\pi}\rightarrow \mathfrak{b} \rightarrow 0$
 be  a short exact sequence,
where $i$ is the inclusion map. For any $b\in B$, choose a $\tilde{b}\in \mathfrak{e}$
such that $\pi(\tilde{b})=b $ in $\mathfrak{b}$. Let $\widetilde{\theta}:k\langle( A\cup B)_+\rangle\rightarrow \mathfrak{e}$
 be the unique homomorphism from $k\langle( A\cup B)_+\rangle$ to $\mathfrak{e}$
such that $\widetilde{\theta}(b)=\tilde{b}$ and $\widetilde{\theta}(a)=a$ for
any $b\in B$ and $a\in A$.

\begin{lemma}\label{le2.4}
Let  $\mathfrak{e}$ be an extension of $\mathfrak{a}$ by
$\mathfrak{b}$,  i.e., $0\rightarrow \mathfrak{a}  \overset{i}\rightarrow \mathfrak{e}\overset{\pi}\rightarrow \mathfrak{b} \rightarrow 0$
is a short exact sequence. Define the left (resp.  right) $k\langle B
_+\rangle$-module operation in $\mathfrak{a}$ as follows:
$$
g\cdot a=\tilde{\theta}(g)a,  \ (resp.  \ a\cdot
g=a\tilde{\theta}(g)),\ a\in  A,\ g\in k\langle B _+\rangle
$$
and
$$
\lfloor \ \rfloor: R\rightarrow  \mathfrak{a},\  \lfloor f\rfloor=
\widetilde{\theta}(f),\ f\in R.
$$
  Then  $\mathfrak{a}$ is a    compatible $k\langle B_+\rangle$-bimodule and $\lfloor \ \rfloor$ is a   factor set of $\mathfrak{b}$ in $\mathfrak{a}$ such that the conditions $(a)$ and $(b)$ in Lemma \ref{le2.3} hold.
 Moreover,
 $$
\mathfrak{e}\cong E_{(\mathfrak{a}, \mathfrak{b},\lfloor\
\rfloor,\cdot)}= k\langle (A\cup B)_+|R_{(\lfloor\  \rfloor
\cdot)}\rangle.
$$

\end{lemma}
{\bf Proof.}
 Since $\mathfrak{a}$ is an ideal of $\mathfrak{e}$, we have
$\tilde{\theta}(g)a\in \mathfrak{a}$, (resp. $a\tilde{\theta}(g)\in \mathfrak{a}$) for any $g\in k\langle B_+\rangle$ and $ a\in  A$.  Thus, the left (resp.  right)  operation is well defined.  It is easy to check that   $\mathfrak{a}$ is a compatible $k\langle B_+\rangle$-bimodule.

 For any $f\in R$, since $\pi(\widetilde{\theta}(f))=f =0$ in $\mathfrak{b}$,  we have
$\widetilde{\theta}(f)\in \mathfrak{a}$. Thus,  $\lfloor \ \rfloor$ is well defined.

$(a)$ \ For any $a\in A, f\in R$, we have
$$
f\cdot a = \tilde{\theta}(f)a=  \{\lfloor f\rfloor a\},\ \ a\cdot f=a\tilde{\theta}(f)=\{ a\lfloor f\rfloor\}.
$$

$(b)$ \ For any intersection  composition $(f_1, f_2)_w$, $f_1, f_2\in R,$ $w=\overline{f_1}c=d\overline{f_2}, c,d\in B^*$  with  $|\overline{f_1}|+|\overline{f_2}|>|w|$, if
$
(f_1,f_2)_w=f_1c-df_2=\sum \alpha_i c_ig_id_i
 $
where each $\alpha_i\in k,  c_i, d_i\in B^*, g_i\in R, c_i\overline{g_i}d_i<w$, then
\begin{eqnarray*}
0&=&\tilde{\theta}(f_1c-df_2-\sum \alpha_i   c_ig_id_i)\\
&=& \tilde{\theta}(f_1)\tilde{\theta}(c)-\tilde{\theta}(d)\tilde{\theta}(f_2)-\sum \alpha_i   \tilde{\theta}(c_i)\tilde{\theta}(g_i)\tilde{\theta}(d_i)\\
&=& \lfloor f_1\rfloor\tilde{\theta}(c)-\tilde{\theta}(d)\lfloor f_2\rfloor-\sum \alpha_i   \tilde{\theta}(c_i) \lfloor g_i\rfloor \tilde{\theta}(d_i)\\
&=& \lfloor f_1\rfloor \cdot c-d\cdot\lfloor f_2\rfloor-\sum \alpha_i  c_i\cdot\lfloor   g_i \rfloor  \cdot d_i.
 \end{eqnarray*}

 Let $\widetilde{\mathfrak{b}}$ be the subspace of $\mathfrak{e}$ spanned by $\tilde{\theta}(Irr(R))=\{\tilde{\theta}(u)|u\in Irr(R)\}$. For any $u, u'\in Irr(R)$, if $
uu'=\sum \delta_iu_i+\sum \beta_j c_jg_jd_j
$
in $k\langle B_+\rangle$, where $\delta_i, \beta_j\in k$,  $u_i\in Irr(R)$, $c_j, d_j\in B^*$ and $g_j\in R$,
then
$$
\tilde{\theta}(u)\tilde{\theta}(u')=\tilde{\theta}(uu')
=\sum \delta_i \tilde{\theta}(u_i)+\sum \beta_j \tilde{\theta}(c_j)\lfloor g_j \rfloor \tilde{\theta}( d_j).
$$
Therefore, $ \mathfrak{e}=\mathfrak{a}\oplus
\widetilde{\mathfrak{b}} $ as vector space with the multiplication:
for any $a,a'\in A$ and $u, u'\in Irr(R)$,
\begin{eqnarray*}
&&(a+\tilde{\theta}(u))(a'+\tilde{\theta}(u'))\\
&=&aa'+ \tilde{\theta}(u)a'
+a\tilde{\theta}(u)+\sum \delta_i \tilde{\theta}(u_i)+\sum \beta_j \tilde{\theta}(c_j)\lfloor g_j \rfloor \tilde{\theta}( d_j)\\
&=&\{aa'\}+u\cdot a' +a\cdot u+\sum \delta_i
\tilde{\theta}(u_i)+\sum \beta_j  c_j\cdot \lfloor g_j \rfloor\cdot
d_j.
\end{eqnarray*}

Define a linear map
$$
\varphi: E_{(\mathfrak{a}, \mathfrak{b},\lfloor\  \rfloor, \cdot )}=\mathfrak{a}\oplus
\mathfrak{b}\rightarrow  \mathfrak{e}=\mathfrak{a}\oplus \widetilde{\mathfrak{b}},\
 a\mapsto a, \  u\mapsto  \tilde{\theta}(u),\ \ a \in A,\ u\in Irr(R).
$$
Therefore,   by Lemma \ref{le2.3}, we can see   that  $\varphi$ is an  algebra isomorphism. \hfill$ \square$\\

\noindent {\bf Remark:} By  Lemma \ref{le2.4},   for a  given  extension  $0\rightarrow \mathfrak{a}  \overset{i}\rightarrow \mathfrak{e}\overset{\pi}\rightarrow \mathfrak{b} \rightarrow 0$, it  is equivalent to some  extension $
0\rightarrow \mathfrak{a} \overset{i_{ \llcorner\lrcorner,\cdot}}{\longrightarrow}  E_{(\mathfrak{a}, \mathfrak{b},\lfloor\
\rfloor,\cdot)} \overset{\pi_{ \llcorner\lrcorner,\cdot}}{\longrightarrow} \mathfrak{b}\rightarrow 0,
$
where  $\cdot$ is a compatible $k\langle B_+\rangle$-bimodule operation on $\mathfrak{a}$ and  $\lfloor\
\rfloor$ is  a factor set such that the conditions  ($a$) and ($b$) in Lemma \ref{le2.3} hold.   If this the case, we call $(\lfloor\  \rfloor, \cdot)$ a special extension pair associated to the extension  $0\rightarrow \mathfrak{a}  \overset{i}\rightarrow \mathfrak{e}\overset{\pi}\rightarrow \mathfrak{b} \rightarrow 0$.  \\

Let  $S$ be  a  Gr\"{o}bner-Shirshov basis  in $k\langle X_+\rangle $. Then for any $f\in k\langle X_+\rangle $, by Theorem \ref{th2.1}, there is a unique $Red_S(f)\in kIrr(S)$ such that  $f=Red_{S}(f)$ in $k\langle X_+|S\rangle$, where $kIrr(S)$ is the $k$-linear space spanned by $Irr(S)$.

Let $\eta: B\rightarrow \mathfrak{a}$ be  a map and
$
\Gamma_{ \eta}:k\langle  (A\cup B)_+\rangle\rightarrow k\langle  (A\cup B)_+\rangle
$
the unique homomorphism  such that $\Gamma_{\eta}(b)=b-\eta(b)$ and $\Gamma_{\eta}(a)=a$ for any $b\in B,\ a\in A$.  If $R_{(\lfloor\  \rfloor', \cdot')}$ is  a Gr\"{o}bner-Shirshov basis in
$k\langle (A\cup B)_+\rangle$, then for any $g\in k\langle  B _+\rangle\subseteq k\langle (A\cup B)_+\rangle$, we have a unique  $Red_{R_{(\lfloor\  \rfloor', \cdot')}} ( g-\Gamma_{ \eta}(g)) \in \mathfrak{a} $.  \\

By   Lemmas \ref{le2.3} and   \ref{le2.4}, we have the following theorem which gives a complete characterization of associative   extensions of    $\mathfrak{a}$ by     $\mathfrak{b}$.

\begin{theorem}\label{th2.5}

Let $\mathfrak{a}$  be  an associative algebra with  a linear basis
$A$ and the    multiplication table: $aa'=\{aa'\}, a, a'\in A$, and
$\mathfrak{b}=k\langle B_+| R\rangle$ an associative algebra
generated by a set $B$ with defining relations $R$.
Assume that $R$  is a minimal Gr\"{o}bner-Shirshov basis in
$k\langle B_+\rangle$ with respect to a monomial order $>$ on $B^*$. Then the following two statements hold.
\begin{enumerate}
\item[(I)]
An associative algebra $\mathfrak{e}$ is an extension of
$\mathfrak{a}$ by $\mathfrak{b}$ if and only if
$$
\mathfrak{e}\cong E_{(\mathfrak{a}, \mathfrak{b},\lfloor\  \rfloor,\cdot)}= k\left\langle
 (A\cup B)_+
 \left|
 \begin{array}{ll}
 aa'-\{aa'\},  & a, a'\in A  \\
 ab-a\cdot b,  ba-b\cdot a,   & a \in A,  b \in B \\
 f- \lfloor f \rfloor, & f \in R
 \end{array}
\right. \right\rangle
$$
for some
compatible $k\langle B_+\rangle$-bimodule  structure in   $\mathfrak{a}$ (denote  the module operation by $\cdot$)  and
some   factor set $ \lfloor\  \rfloor : R \rightarrow \mathfrak{a} $
of $\mathfrak{b}$ in $\mathfrak{a}$ such that
\begin{enumerate}
\item[(a)]
 For any $a\in A, f\in R$,
 $
f\cdot a =\{\lfloor f\rfloor a\},\ \ a\cdot f=\{ a\lfloor f\rfloor\}
$
in   $\mathfrak{a}$.

\item[(b)] For any intersection  composition $(f_1, f_2)_w$, $f_1, f_2\in R,$
$w=\overline{f_1}c=d\overline{f_2},\ c,d\in B^*$  with  $|\overline{f_1}|+|\overline{f_2}|>|w|$, if
$
(f_1,f_2)_w=f_1c-df_2=\sum \alpha_i c_ig_id_i
 $
where each $\alpha_i\in k,  c_i, d_i\in B^*, g_i\in R, c_i\overline{g_i}d_i<w$, then in $\mathfrak{a}$,
$$
\lfloor f_1\rfloor \cdot c-d\cdot\lfloor f_2\rfloor=\sum \alpha_i
c_i\cdot\lfloor   g_i \rfloor  \cdot d_i.
$$
\end{enumerate}

\item[(II)] Two
extensions $\mathfrak{e}$ and $\mathfrak{e'} $ of $\mathfrak{a}$ by $\mathfrak{b} $ are   equivalent  if and  only if  there exists a   map  $\eta:B\rightarrow
\mathfrak{a}$  such that    in $\mathfrak{a}$
\begin{eqnarray}
&&  a\cdot' b-a\cdot b=\{a\eta(b)\},  \
   b\cdot' a-b\cdot a= \{ \eta(b)a\} \label{c3}\\
 &&\lfloor f\rfloor' -\lfloor f\rfloor = Red_{R_{(\lfloor\  \rfloor', \cdot')}} ( f-\Gamma_{ \eta}(f)) \label{c5}
\end{eqnarray}
for any $ a\in A,\ b\in B $ and    $f\in R$,  where $(\lfloor\ \rfloor,\cdot)\ (resp.  \ (\lfloor\ \rfloor', \cdot'))$ is a special extension pair  associated to $e$ (resp. $e'$).
\end{enumerate}

\end{theorem}
{\bf Proof.}   (I) This part follows from  Lemmas \ref{le2.3} and \ref{le2.4}.

(II) ($\Longrightarrow$)  Assume that $\varphi$ is a homomorphism and the following diagram is commutative:
\begin{eqnarray}\label{tu6}
\xymatrix{
 0 \ar[r]& \mathfrak{a} \ar[r]^{i_{\llcorner\lrcorner, \cdot}\ \ \ \ } \ar@{=}[d]_{}& E_{(\mathfrak{a}, \mathfrak{b},
 \lfloor\  \rfloor,\cdot)}\ar[r]^{\ \ \ \ \ \  \pi_{\llcorner\lrcorner, \cdot}}\ar[d]^{\varphi} &
 \mathfrak{b} \ar@{=}[d]^{} \ar[r]&0    \\
        0 \ar[r] & \mathfrak{a} \ar[r]^{i_{\llcorner\lrcorner', \cdot'}\ \ \ \ }  & E_{(\mathfrak{a},
       \mathfrak{b},\lfloor\  \rfloor', \cdot')}\ar[r]^{\ \   \ \ \ \   \pi_{\llcorner\lrcorner', \cdot'}}&
       \mathfrak{b}  \ar[r]&0 }
\end{eqnarray}
For any $b\in B$, by the  commutative diagram, we have
$$
 \pi_{\llcorner\lrcorner', \cdot'} (\varphi(b))= \pi_{\llcorner\lrcorner, \cdot}(b)=b= \pi_{\llcorner\lrcorner', \cdot'}(b).
$$
Thus, $b-\varphi(b)\in \mathfrak{a}$. Define $\eta: B\rightarrow
\mathfrak{a},\ b\mapsto b-\varphi(b)$. Then, it is easy to check that $(\ref{c3})$ and $(\ref{c5})$ are true.

($\Longleftarrow$) Let
$
\widetilde{\varphi}:k\langle(  A\cup B)_+\rangle  \rightarrow  E_{(\mathfrak{a}, \mathfrak{b},\lfloor\  \rfloor', \cdot')}
$
be the unique homomorphism from $k\langle(  A\cup B)_+\rangle $  to
$E_{(\mathfrak{a}, \mathfrak{b},\lfloor\  \rfloor', \cdot')}$ such
that $\widetilde{\varphi}(a)=a$ and $\widetilde{\varphi}(b)=b-\eta(b)$   for any $a\in A, b\in
B$. By (\ref{c3}) and (\ref{c5}), we have  $\widetilde{\varphi}(h)=0$ for any $h\in R_{(\lfloor\ \rfloor,\cdot)}$.   Therefore, $\widetilde{\varphi}$ induces a  homomorphism
$
\varphi:E_{(\mathfrak{a}, \mathfrak{b},\lfloor\  \rfloor, \cdot)}\rightarrow
E_{(\mathfrak{a}, \mathfrak{b},\lfloor\  \rfloor', \cdot')}
$
which makes the diagram
(\ref{tu6}) commutative.
\hfill $\square$\\

\subsection{Some special cases of extensions of associative algebras}

Let $\mathfrak{b}$ be an associative algebra with a well-ordered linear basis $B$ and the multiplication table:  $bb'=\{bb'\}, b,b'\in B$.
  It is easy to see that $R=\{bb'-\{bb'\}|b,b'\in B\}$ is a minimal Gr\"{o}bner-Shirshov basis with respect to the deg-lex order on $ B ^*$.

For any $b,b'\in B$, denote
$
(b,b')=\lfloor bb'-\{bb'\}\rfloor,
$
which can be extended to be a bilinear map
$
( , ):\mathfrak{b}\times \mathfrak{b}\rightarrow \mathfrak{a}.
$

 By  Theorem \ref{th2.5}, we have the following corollary.

 \begin{corollary} \label{th2.7} Let $\mathfrak{a}$ (resp. $\mathfrak{b}$)  be  an
 associative algebra with  a linear basis $A$ (resp. $B$)  and the
 multiplication table: $aa'=\{aa'\}, a, a'\in A$ (resp. $bb'=\{bb'\}, b,b'\in B$).
 Then an associative  algebra $\mathfrak{e}$ is an extension of $\mathfrak{a}$ by $\mathfrak{b}$ if and only if
  $$
\mathfrak{e}\cong  k\left\langle
 (A\cup B)_+
 \left|
 \begin{array}{ll}
 aa'-\{aa'\},  & a, a'\in A  \\
 ab-a\cdot b,  \ ba-b\cdot a,  & a \in A,  b \in B \\
 bb'- \{bb'\}-(b,b'),  & b, b' \in B
 \end{array}
\right. \right\rangle
$$
for some  compatible $k\langle  B _+\rangle$-bimodule  structure in  $\mathfrak{a}$ (denote  the module operation by $\cdot$)
 and  some    bilinear map
$
( , ):\mathfrak{b}\times \mathfrak{b}\rightarrow \mathfrak{a}
$
such that  for any $a\in A, b,b', b''\in B$,

$$
(bb'-\{bb'\})\cdot a =\{(b,b')a\}, \ \ a\cdot (bb'-\{bb'\})=\{a(b,b')\},
$$

$$
(\{bb'\},b'')+ ( b,b')\cdot b''-(b, \{b'b''\})-b\cdot ( b',b'')=0.
$$

\end{corollary}
\ \

Now, we consider another special case of Theorem \ref{th2.5}:  $\mathfrak{a}^2=0$.

\begin{corollary}\label{th2.8}(\cite{chen09}) Let $\mathfrak{a}$  be  an associative algebra with a linear basis $A$ and
  $\mathfrak{b}=k\langle B_+| R\rangle$, where  $R$  is a  minimal Gr\"{o}bner-Shirshov basis in $k\langle B_+\rangle$.
If
$\mathfrak{a}^2=0$, then an associative  algebra $\mathfrak{e}$ is an extension of $\mathfrak{a}$ by $\mathfrak{b}$ if and only if
  $$
\mathfrak{e}\cong E_{(\mathfrak{a}, \mathfrak{b},\lfloor\  \rfloor,\bullet)}= k\left\langle
 (A\cup B)_+
 \left|
 \begin{array}{ll}
 aa',  & a, a'\in A  \\
 ab-a\bullet b,  \ ba-b\bullet a,  & a \in A,  b \in B \\
f-\lfloor f\rfloor,  & f \in R
 \end{array}
\right. \right\rangle
$$
for some
  $\mathfrak{b}$-bimodule structure in $\mathfrak{a}$  (denote the module operation by $\bullet$) and    some   factor set
  $\lfloor\ \rfloor: R\rightarrow \mathfrak{a}$  of $\mathfrak{b}$ in $\mathfrak{a}$ such that
\begin{enumerate}
\item[(b)] For any intersection  composition $(f_1, f_2)_w$, $f_1, f_2\in R,$ $w=\overline{f_1}c=d\overline{f_2}, c,d\in B^*$  with  $|\overline{f_1}|+|\overline{f_2}|>|w|$, if
$
(f_1,f_2)_w=f_1c-df_2=\sum \alpha_i c_ig_id_i
 $
where each $\alpha_i\in k,  c_i, d_i\in B^*, g_i\in R, c_i\overline{g_i}d_i<w$, then
$$
 \lfloor f_1\rfloor \bullet c-d\bullet\lfloor f_2\rfloor=\sum \alpha_i  c_i\bullet \lfloor   g_i \rfloor  \bullet  d_i.
$$
\end{enumerate}

\end{corollary}
{\bf Proof.}  $(\Longleftarrow)$ If $\mathfrak{a}^2=0$  and $\mathfrak{a}$  is a $\mathfrak{b}$-bimodule with the module operation $\bullet$,
then it is easy to see that  $\mathfrak{a}$  is a compatible $k\langle B_+\rangle$-bimodule with the module  operation
$$
a\cdot g=a\bullet (g+Id(R)), \ g\cdot a =(g+Id(R))\bullet a, \ a\in A, g\in k\langle B_+\rangle.
$$
For any $a\in A, f\in R$, since $\mathfrak{a}^2=0$, we have
$$
f\cdot a =(f+Id(R))\bullet a =0=\{\lfloor f\rfloor a\},\ \ a\cdot f=a\bullet (f+Id(R))=0=\{ a\lfloor f\rfloor\}.
$$
Therefore,  by Theorem \ref{th2.5}, $\mathfrak{e}$ is an extension of $\mathfrak{a}$ by $\mathfrak{b}$.

$(\Longrightarrow)$  If $\mathfrak{e}$ is an extension of $\mathfrak{a}$ by $\mathfrak{b}$, then by Theorem \ref{th2.5}, $\mathfrak{a}$  is a compatible $k\langle B_+\rangle$-bimodule and $ \mathfrak{e}\cong E_{(\mathfrak{a}, \mathfrak{b},\lfloor\  \rfloor,\cdot)}$ for some factor set $\lfloor\ \rfloor: R\rightarrow  \mathfrak{a}$ such that the conditions $(a)$ and $(b)$ in Theorem  \ref{th2.5} hold.

Since $\mathfrak{a}^2=0$,     for any $f\in R, a\in A$, we have
$$
f\cdot a=\{\lfloor f\rfloor a \}=0, \  a\cdot f = \{ a \lfloor f  \rfloor \}=0.
$$
Thus,  $\mathfrak{a}$ is  $\mathfrak{b}$-bimodule   with the operation:
$$
a\bullet (g+Id(R))= a\cdot g, \ (g+Id(R))\bullet a= g\cdot a, \ a\in A, g\in k\langle B_+\rangle.
$$
It is easy to see that  $ \mathfrak{e}\cong E_{(\mathfrak{a}, \mathfrak{b},\lfloor\  \rfloor,\cdot)}=E_{(\mathfrak{a}, \mathfrak{b},\lfloor\  \rfloor,\bullet)}$ and the condition $(b)$ holds.\hfill $\square$\\

\ \

By Corollaries   \ref{th2.7} and \ref{th2.8}, we have the following corollary.

\begin{corollary}(\cite{lan}) Let $\mathfrak{a}$ be an associative algebra with  a linear basis $A$  and  $\mathfrak{b}$  be  an
 associative algebra with  a linear basis $B$   and the
 multiplication table:   $bb'=\{bb'\}, b,b'\in B$.  If   $\mathfrak{a}^2=0$, then
an algebra $\mathfrak{e}$  is an extension  of $\mathfrak{a}$ by $\mathfrak{b}$  if and only if
$$
\mathfrak{e}\cong  k\left\langle
 (A\cup B)_+
 \left|
 \begin{array}{ll}
 aa',  & a, a'\in A  \\
 ab-a\bullet b,  \ ba-b\bullet a,  & a \in A,  b \in B \\
bb'-\{bb'\}-(b,b'),   & b, b' \in B
 \end{array}
\right. \right\rangle
$$
for some $\mathfrak{b}$-bimodule structure in $\mathfrak{a}$  (denote the module operation by $\bullet$) and   some  bilinear map $( , ):\mathfrak{b}\times \mathfrak{b}\rightarrow \mathfrak{a}$ such that for any $b,b',b''\in B$,
$$
( \{b b'\}, b'')+(b,b')\bullet b''-(b, \{b' b''\})-b\bullet(b',b'' )=0.
$$
 \end{corollary}

\subsection{Applications}

 In this subsection, we  assume that  $\mathfrak{a}$ is  an  associative algebra with a well-ordered linear basis   $A$ and  multiplication table:   $ aa'= \{aa'\}, a,a'\in A$. Let $\mathfrak{b}=k\langle B_+| R\rangle$ with $B$ a well-ordered set. We assume  that the order on $(A\cup B)^*$ is the deg-lex order  such that $b>a$ for any $b\in B, a\in A$.

By Theorem \ref{th2.5} and direct computations, we have the following examples.
\begin{example} Let   $\mathfrak{b}=k\langle \{b\}_+| b^n=f(b)\rangle$ be a cyclic associative algebra, where $n$  is a natural number and $f(b)$ is a polynomial of degree less than $n$ such that $f(0)=0$.
Then an associative algebra $\mathfrak{e}$ is   an extension of   $\mathfrak{a}$ by $\mathfrak{b}$ if and only if
$$
\mathfrak{e}\cong k\langle (A\cup \{b\})_+|aa'- \{aa'\}, b^n-f(b)-\tilde{a}, ba-b\cdot a, ab-a\cdot b, a,a'\in A\rangle
$$
for some
compatible $k\langle  \{b\}_+\rangle$-bimodule  structure in   $\mathfrak{a}$   and some  $\tilde{a}\in \mathfrak{a}$ such that for  any $a \in A$,
$$
b^n\cdot a-f(b)\cdot a=\{\tilde{a}a \}, a\cdot b^n-a\cdot f(b)=\{a\tilde{a} \}, \  b^i\cdot \tilde{a}= \tilde{a}\cdot b^i, 1\leq i\leq n-1.
$$
 \end{example}

\begin{example}Let  $\mathfrak{b}=k\langle B_+| bb'-b'b, b>b', b, b'\in B\rangle$ be the   free commutative algebra on $B$. Then an associative algebra $\mathfrak{e}$ is
 an extension of   $\mathfrak{a}$ by $\mathfrak{b}$ if and only if
$$
\mathfrak{e}\cong k\langle (A\cup B)_+|aa'-\{aa'\}, ba-b\cdot a, ab-a\cdot b, bb'-b'b-(b,b'), a,a'\in A, b, b'\in B\rangle
$$
for some  compatible $k\langle  B_+\rangle$-bimodule  structure in $\mathfrak{a}$  and
  some factor set $\{(b,b')\in \mathfrak{a} |b,b'\in B, b>b'\}$ such that  for any $a\in A$ and $b,b', b''\in B, b>b'>b''$,
$$
b\cdot( b'\cdot a)-b'\cdot (b\cdot a)=\{(b,b')a\},\ (a \cdot  b )\cdot b'-(a\cdot  b') \cdot b=\{a(b,b')\},
$$
$$
(b',b'')\cdot b +b'\cdot(b,b'')+(b,b')\cdot b''-(b,b'')\cdot b'-b''\cdot (b,b')-b\cdot (b',b'')=0.
$$
 \end{example}

\section{Extensions of Lie algebras}

\subsection{Composition-Diamond  lemma  for   Lie algebras}
In  this subsection, we review the Gr\"obner-Shirshov bases   theory for Lie algebras,  which can be found in \cite{b07,BC14, Sh, Shir3}.

Let   $X$ be a set which is equipped with  a well order $>_X$ and    $\succ$ (resp. $>$) be the lex order (resp. deg-lex) order on $X^*$ with respect to the order $>_X$.   Let $X^{**}$ be the  set of all  non-associative words on  $X$.

With the lex order $\succ$ on $X^*$, we   give the definition of   Lyndon-Shirshov words, see \cite{b07,BC14, Sh, Shir3} for details.

A nonempty  word $w\in X^*$ is called an associative Lyndon-Shirshov word on $X$ if $ w=uv, u, v\neq 1\Rightarrow w \succ vu$.
 A non-associative word $(u)\in X^{**}$ is a non-associative Lyndon-Shirshov word, denoted by  $[u]$, if
  \begin{itemize}
    \item [(a)] $u$ is an associative Lyndon-Shirshov word;
    \item [(b)] if $(u)=((u_1)(u_2))$, then both $(u_1)$ and $(u_2)$ are non-associative Lyndon-Shirshov words;
    \item [(c)] if $(u_1)=((u_1')(u_1'')) $, then $u_1''\preceq u_2$.
  \end{itemize}
Denote $ALSW(X)$ (resp. $NLSW(X)$) to be  the set of all  associative (resp.  non-associative) Lyndon-Shirshov words on $X$.

In fact, for any associative Lyndon-Shirshov word $w$, there is a unique bracketing way $[w]$ such that $[w]$ is a non-associative Lyndon-Shirshov word.
Let $k\langle X\rangle$   be the free associative algebra with identity  on $X$ over a field $k$ and $Lie(X)$ the Lie subalgebra of
$k\langle X\rangle$ generated by $X$ under the Lie bracket $(uv)=uv-vu$. It is well known that $Lie(X)$ is a free Lie algebra on $X$ with a linear basis  $NLSW(X)$.

Let $(u) = (x_1x_2\cdots x_m) \in X^{**}$, where each $x_i\in X$. Denote $cont( u )$ the multiset $\{x_1, x_2, \ldots, x_m\}$.

\begin{lemma}\label{le2.1} (\cite{b07,Sh, Shir3})
For any non-associative word $(u)\in X^{**}$, $(u)$ has a unique representation
$$
(u)=\alpha_1[u_1]+\alpha_2[u_2]+\cdots + \alpha_n [u_n],
$$
where each $\alpha_i\in k$ and  $u_i\in ALSW(X)$, $u_1>u_2>\cdots u_n$ and $cont( u )=cont(u_i)$.
\end{lemma}

For any $f\in k\langle X\rangle$,    denote $\overline{f}$ the leading word of $f$ with respect to the deg-lex order $>$ on $X^*$. If  $[w]\in NLSW(X)$, then we have $\overline{[w]}=w$.

If $w=avb$ with $w$ and $v$  associative Lyndon-Shirshov words, then
$$
[w]=[a[vc]d],
$$
where $b=cd$ for some $c,d\in X^*$.  Denote
$$
[w]_{v}=[w]|_{[vc]\mapsto [\cdots[[[v][c_1]][c_2]]\cdots [c_n]]},
$$
where $c=c_1c_2\cdots c_n$, each  $c_i$ is an  associative Lyndon-Shirshov word and $c_1\preceq c_2\preceq \cdots\preceq c_n$. Then $\overline{[w]_v}=w$.

Let $s\in Lie(X)$ be a monic polynomial.  If $a\overline{s}b$ is an associative Lyndon-Shirshov word, define
$$
[asb]_{\overline{s}}=[a\overline{s}b]_{\overline{s}}|_{[\overline{s}]\mapsto s}
$$
which is called a  special normal $s$-word.

Let $f,  g \in Lie(X)$ be  monic. There are two kinds of compositions.

(i)  If  $w=\bar{f}a=b\bar{g}$ for some $a,b\in X^*$  with  $|\bar{f}|+|\bar{g}|>|w|$, then
$
\langle f,g\rangle_w=[fa]_{\overline{f}} - [bg]_{\overline{g}}
$
 is called the intersection composition of $f$ and $g$ with respect to the  ambiguity $w$.

(ii)  If $w=\overline{f}=a\overline{g}b$ for some  $a,b\in X^*$, then
$
\langle f,g\rangle_w=f-[agb]_{\overline{g}}
$
is called the inclusion composition of $f$ and $g$ with respect to  the ambiguity $w$.

 If $S$ is a monic subset of $Lie(X)$, then
the composition $\langle f,g\rangle_w$ is called trivial modulo $(S, w)$,  if
$\langle f,g\rangle_w=\sum\alpha_i [a_i s_i b_i]_{\overline{s_i}}$, where each $\alpha_i\in k$,
$a_i,b_i\in X^{*}, \ s_i\in S$, $a_i \overline{s_i} b_i\in ALSW(X)$ and $a_i \overline{ s_i }b_i<w$. If
this is the case, then we write
$$
\langle f,g\rangle_w\equiv0\  mod(S,w).
$$

A monic set $S$ is called a
Gr\"{o}bner-Shirshov basis  in   $Lie(X)$ if any
composition $\langle f,g\rangle_w$ of $f, g\in S$ is trivial modulo $(S,w)$.

A Gr\"{o}bner-Shirshov basis $S$ in   $Lie(X)$ is minimal if there is no inclusion composition in $S$. Let $\mathfrak{c}=Lie( X | R)$  be a Lie algebra
generated by  $X$  and defining relations $R$.
 Then, one can find a minimal
Gr\"{o}bner-Shirshov basis $S$ in   $Lie(X)$ such that $\mathfrak{c}=Lie( X | R)=Lie( X |S)$, see, \cite{BC14}.

\begin{theorem} (\cite{b07,Sh, Shir3})\label{CD_lemmaLie}
(Composition-Diamond lemma for Lie algebras) Let $S\subseteq Lie(X)$ be a non-empty monic set, $Id(S)$    the ideal of $Lie(X)$ generated by $S$ and $>$ the deg-lex order on $X^*$.
 Then the following statements are equivalent:
\begin{enumerate}
\item[(i)] $S $ is a Gr\"{o}bner-Shirshov basis  in $Lie(X)$.
\item[(ii)] $f\in Id(S)\Rightarrow \bar{f}=a\bar{s}b \in ALSW(X)$
for some $s\in S$ and $a,b\in  X^*$.
\item[(iii)] The set
$$
Irr(S)= \{[w]|w \in ALSW(X), w\neq a\bar{s}b, s\in S, a,b \in X^* \}
$$
is a  linear basis of the Lie  algebra $Lie(X|S)=Lie(X)/Id(S)$.
\end{enumerate}
\end{theorem}

\subsection{A characterization of extensions of Lie algebras}

Let $\mathfrak{a},\mathfrak{b}$ be   Lie  algebras over a  field $k$.
In this subsection, by using the the Gr\"{o}bner-Shirshov bases theory   for Lie algebras,  we give  a complete characterization of extensions of      $\mathfrak{a}$ by    $\mathfrak{b}$,  where   $\mathfrak{b}$ is presented by generators and relations.

\begin{definition}
Let $\mathfrak{c}$ and  $\mathfrak{d}$ be two Lie algebras. Then $\mathfrak{c}$ is called a  compatible $\mathfrak{d}$-module if  $\mathfrak{c}$ is   a  $\mathfrak{d}$-module and for any $c, c'\in \mathfrak{c}$ and $d\in \mathfrak{d}$,
$$
d\cdot [cc']=[(d\cdot c) c']+[c(d\cdot c')],
$$
where $\cdot$ is the $\mathfrak{d}$-module operation.

\end{definition}
Let   $\mathfrak{c}$ and  $\mathfrak{d}$ be two Lie algebras and $\mathfrak{c}$    a  $\mathfrak{d}$-module relating to  the representation
$
\phi:  \mathfrak{d} \rightarrow  \mathfrak{gl(c)},
$
where $ \mathfrak{gl(c)}$ is the Lie algebra   of all linear transformations of     $\mathfrak{c}$.
Then  $\mathfrak{c}$ is  a  compatible $\mathfrak{d}$-module if and only if $\phi(d)$ is a derivation of  $\mathfrak{c}$ for any $d\in \mathfrak{d}$.
\\

Let $\mathfrak{a}$ be  a Lie algebra with a linear basis $A $ and the multiplication table: $[aa']-[aa']_\mathfrak{a}, a, a'\in A $ and  $\mathfrak{b}=Lie(B|R)$ a Lie algebra generated by  generators  $B$  and   defining  relations $R$.

Let  $>$ be a well order on $B$    and  extend it  to  a deg-lex order $>$ on $(A\cup B)^*$ such that $b> a$  for any $b\in B, a\in A$.
For convenience, we may assume that $R$  is a  minimal Gr\"{o}bner-Shirshov basis in $Lie(B)$ with respect to the deg-lex order on $B^*$.
 For any $f\in R$, denote
$
r_{_f}=f-[\overline{f}],
$
where $\bar{f}$ is the   leading word of $f$ with respect to the deg-lex order $>$.
  Then for any $f\in Lie(B)$, by Theorem \ref{CD_lemmaLie}, there is a unique $Red_R(f)\in kIrr(R)$ such that  $f=Red_{R}(f)$ in $Lie(B|R)$, where $kIrr(R)$ is the $k$-linear space spanned by $Irr(R)$. Note that for any $f,g\in Lie(B)$, we have
$$
Red_{R}((fg))=Red_{R}((Red_{R}(f)g))=Red_{R}((Red_{R}(f)Red_{R}(g))).
$$

Let   $\lfloor\  \rfloor : R \rightarrow \mathfrak{a}$  be a map, which is called a factor set of $\mathfrak{b}$ in $\mathfrak{a}$. Assume that
$\mathfrak{a}$ is a compatible $Lie(B)$-module with the module operation $\cdot$.

It is easy to see that $S_{(\cdot)}=\{[ba]-b\cdot a |a\in A, b\in B\}$ is a Gr\"{o}bner-Shirshov basis in $Lie(A\cup B)$.
If  $g= \sum  \alpha_i (v_ia_iw_i)$, where each $\alpha_i\in k$ and $(v_ia_iw_i)\in (A\cup B)^{**},  v_i, w_i \in B^*, a_i\in A$, then  we can see that
$Red_{S_{(\cdot)}}(g)\in  \mathfrak{a}$.

  Define
$$
LieE_{(\mathfrak{a}, \mathfrak{b},\lfloor\  \rfloor, \cdot)}=Lie(A\cup B|LieR_{(\lfloor\  \rfloor, \cdot)}),
$$
where  $LieR_{(\lfloor\  \rfloor, \cdot)}$ is
the set consisting of the following   polynomials in $Lie(A\cup B)$:
\begin{eqnarray*}
f_{aa'}:&=&[aa']-[aa']_\mathfrak{a},\ \  a, a'\in A, a>a',\\
t_{ba}:&=&[ba]-b\cdot a,\ \ \ \ \ b\in B, a\in A,\\
h_f:&=&f-\lfloor f\rfloor, \ \ \ \ \ \ \ \ \ f\in R.
\end{eqnarray*}

\begin{lemma}\label{le3.3} If  $((v)a)\in (A\cup B)^{**},  v  \in B^*, a\in A$, then
$
Red_{S_{(\cdot)}}(((v)a))= (v)\cdot a.
$
It follows that for any $f\in R,\ a\in A$,
$
Red_{S_{(\cdot)}}((fa)) = f\cdot a.
$
\end{lemma}
{\bf Proof.} Induction on $|v|$. If $|v|=0$ or 1, then the result is trivial. Assume that the result is true for $(v)$ with   $|v|<n$.  Let $|v|=n\geq2$ and $(v)=((v_1)(v_2))$. Since
$$
((v)a) = (((v_1)(v_2))a)
 =   (((v_1)a)(v_2)) +  ((v_1) ((v_2)a))\\
$$
and by induction
$$
Red_{S_{(\cdot)}}(((v_1)a)))=(v_1)\cdot a, \ Red_{S_{(\cdot)}}(((v_2)a))=(v_2)\cdot a,
$$
we have
\begin{eqnarray*}
&&Red_{S_{(\cdot)}}(((v)a)) \\
 &=&
 Red_{S_{(\cdot)}} ( (((v_1)a)(v_2)) +  ((v_1) ((v_2)a)))\\
 &=&  Red_{S_{(\cdot)}} ( ( Red_{S_{(\cdot)}} (((v_1)\cdot a))(v_2))+  ((v_1)Red_{S_{(\cdot)}}(((v_2)\cdot a))))\\
&= &- (v_2)\cdot((v_1)\cdot a)+ (v_1)\cdot ((v_2)\cdot a)\\
&=& ((v_1)(v_2))\cdot (a)
=(v)\cdot a.
 \end{eqnarray*}
Therefore,
$
Red_{S_{(\cdot)}}(((v)a))= (v)\cdot a.
$
\hfill $\square$

 \begin{lemma}\label{le3.5}
 Let $u\in ALSW(B)$ and $a\in A$. Then in $Lie(A\cup B)$,
 $$
 ([u]a)=[ua]+\sum \alpha_i[v_iaw_i],
 $$
 where $\alpha_i\in k$, $v_i, w_i\in B^*$ and $v_iaw_i<ua$.
 \end{lemma}
 {\bf Proof.}
By Lemma \ref{le2.1}, we have
$$
 ([u]a)=\sum_{i=1}^n   \alpha_i[v_iaw_i],
$$
where each $\alpha_i\in k$, $v_i, w_i\in B^*$ and  $v_1aw_1> v_2aw_2>\cdots> v_naw_n$. Since  $([u]a)=[u]a-a[u]$ and $\overline{[u]}=u$, we have  $\overline{([u]a)}=ua$ and the coefficient of  $\overline{([u]a)}$ is $1$. Thus  $v_1aw_1=ua$ and $\alpha_1=1$. Therefore,
$
 ([u]a)=[ua]+\sum \alpha_i[v_iaw_i].
 $
\hfill $ \square$
 \begin{lemma} \label{th3.5} With the above notations,
$LieR_{(\lfloor\  \rfloor, \cdot)}$ is a Gr\"{o}bner-Shirshov basis  in $Lie(A\cup B)$ if and only if
\begin{enumerate}
\item[(a)]  For any $a\in A,  f\in R$,   $f\cdot a=[\lfloor f\rfloor a]_{\mathfrak{a}}$   in  $\mathfrak{a}$.
\item[(b)] For any intersection  composition $\langle f_1, f_2 \rangle_w$, where $ f_1, f_2\in R,  w=\overline{f_1}c=d\overline{f_2}$, $c, d \in B^* $ and $|\overline{f_1}|+|\overline{f_2}|>|w|$,   if
$$
\langle f_1,f_2\rangle_w=[f_1c]_{\overline{f_1}}-[df_2]_{\overline{f_2}}=\sum\alpha_i[c_ig_id_i]_{\overline{g_i}},
$$
where each $\alpha_i\in k$, $c_i, d_i\in B^*$, $g_i\in R$ and  $c_i\overline{g_i}d_i<w$,
then in  $\mathfrak{a}$,
$$
 Red_{S_{(\cdot)}}([\lfloor f_1 \rfloor c]_{\overline{f_1}} - [d \lfloor f_2\rfloor ]_{\overline{f_2}})=\sum\alpha_i Red_{S_{(\cdot)}}([c_i  \lfloor g_i\rfloor   d_i]_{\overline{g_i}}).
$$

\end{enumerate}
  Moreover, if this is the case,
$
  LieE_{(\mathfrak{a}, \mathfrak{b},\lfloor\  \rfloor,\cdot)}
$
is an extension of $\mathfrak{a}$ by $\mathfrak{b}$.
\end{lemma}
{\bf Proof.} Since $R$ is a minimal Gr\"{o}bner-Shirshov basis in $Lie(B)$,
the only  possible compositions in $LieR_{(\lfloor\  \rfloor, \cdot)}$ are related to the following ambiguities:
$$
aa'a'',\ baa',\  \overline{f}a=vba,\ \overline{f_1}c=d\overline{f_2},
$$
where  $f_1,f_2, f  \in R$, $v,c, d\in B^*, b\in B$,  $a,a', a''\in A, a>a'>a''$ and $|\overline{f_1}|+|\overline{f_2}|>|\overline{f_1}c|$.\\

If $w=aa'a''$, then $\langle f_{aa'},f_{a'a''}\rangle_{w}
\equiv 0\ \ mod (LieR_{(\lfloor\  \rfloor, \cdot)},w)$    by the  Jacobi   identity.

If $w=baa'$,  then since  $\mathfrak{a}$ is a compatible $Lie(B)$-module,  we have
\begin{eqnarray*}
\langle t_{ba},f_{aa'}\rangle_{w}&=&[t_{ba}a']_{\overline{t_{ba}}}-[bf_{aa'}]_{\overline{f_{aa'}}}\\
&=& (([ba]-b\cdot a)a')-(b([aa']-[aa']_\mathfrak{a}))\\
&=&   ([ba]a')-((b\cdot a)a')-[b[aa']]+[b[aa']_\mathfrak{a}]\\
&\equiv& [[ba']a]-[(b\cdot a) a']_\mathfrak{a}+[b[aa']_\mathfrak{a}]\\
&\equiv&  -[a(b\cdot a')]_\mathfrak{a}-[(b\cdot a) a']_\mathfrak{a}+b\cdot [aa']_\mathfrak{a}\\
&\equiv& 0\ \ mod(LieR_{(\lfloor\  \rfloor, \cdot)},w).
\end{eqnarray*}

If  $w=\overline{f}a=vba$, then   by Lemma \ref{le3.5},
$$
([\bar{f}]a)=([vb]a)=[vba]+\sum\alpha_i[v_iaw_i],
$$
where each $\alpha_i\in k,\  v_i, w_i\in  B^*$, and $ v_iaw_i <w$. Thus,
by  Lemma \ref{le3.3},
\begin{eqnarray*}
[\overline{f}]\cdot a= [vb]\cdot a =Red_{S_{(\cdot)}}(([\bar{f}]a))= Red_{S_{(\cdot)}} ([vba])+\sum\alpha_iRed_{S_{(\cdot)}} ([v_i a  w_i]).
\end{eqnarray*}
Therefore, by Lemma \ref{le3.3},
\begin{eqnarray*}
 \langle h_{f},t_{ba}\rangle_w
&=& [(f-\lfloor f\rfloor)a]_{\overline{f}}-[v([ba]-b\cdot a)]_{\overline{t_{ba}}}\\
&=&([\overline{f}]a)+(r_{_f}a)-(\lfloor f\rfloor a)-[v[ba]]+[v ( b\cdot a)]_{_{\overline{t_{ba}}}} \\
&\equiv&\sum\alpha_iRed_{S_{(\cdot)}}  ([v_i a  w_i])+ r_{_f}\cdot a -[\lfloor f\rfloor a]_{\mathfrak{a}}+ Red_{S_{(\cdot)}} ([vba]) \\
&\equiv&[vb]\cdot a +r_{_f}\cdot a -[\lfloor f\rfloor a]_{\mathfrak{a}} \\
&\equiv& f\cdot a -[\lfloor f\rfloor a]_{\mathfrak{a}}\ \ mod(LieR_{(\lfloor\  \rfloor, \cdot)}, w).
\end{eqnarray*}
Then  $\langle h_{f},t_{ba}\rangle_w\equiv0\  mod(LieR_{(\lfloor\  \rfloor, \cdot)},w)$ if and only if
$
f\cdot a=[\lfloor f\rfloor a]_{\mathfrak{a}}.
$

If $w=\overline{f_1}c=d\overline{f_2}$, $f_1, f_2\in R, c, d \in B^*$ and $|\overline{f_1}|+|\overline{f_2}|>|w|$, then
\begin{eqnarray*}
\langle h_{f_1},h_{f_2}\rangle_w
&=&[(f_1-\lfloor f_1\rfloor)c]_{\overline{f_1}}-[d(f_2-\lfloor f_2\rfloor)]_{\overline{f_2}}\\
&=& [f_1c]_{\overline{f_1}}-[df_2]_{\overline{f_2}}+([d\lfloor f_2\rfloor]_{\overline{f_2}}-[\lfloor f_1\rfloor c]_{\overline{f_1}}).
\end{eqnarray*}
Since  $R$ is a Gr\"{o}bner-Shirshov basis in $Lie(B)$, we have
$$
\langle f_1,f_2\rangle_w=[f_1c]_{\overline{f_1}}-[df_2]_{\overline{f_2}}=\sum\alpha_i[c_ig_id_i]_{\overline{g_i}},
$$
where each $\alpha_i\in k$,  $c_i, d_i\in B^*$,  $g_i\in R$  and  $c_i\overline{g_i}d_i<w$.
Therefore,
\begin{eqnarray*}
&& \langle h_{f_1},h_{f_2}\rangle_w\\
&=&\sum\alpha_i[c_i(g_i-\lfloor g_i\rfloor) d_i]_{\overline{g_i}}+\sum\alpha_i[c_i\lfloor g_i\rfloor d_i]_{\overline{g_i}}+[d\lfloor f_2\rfloor]_{\overline{f_2}}-[\lfloor f_1\rfloor)c]_{\overline{f_1}}\\
&\equiv& \sum\alpha_i[c_i\lfloor g_i\rfloor d_i]_{\overline{g_i}}+[d\lfloor f_2\rfloor]_{\overline{f_2}}-[\lfloor f_1\rfloor c]_{\overline{f_1}} \\
&\equiv&\sum\alpha_i Red_{S_{(\cdot)}} ([c_i\lfloor g_i\rfloor d_i]_{\overline{g_i}})- Red_{S_{(\cdot)}} ( [\lfloor f_1\rfloor c]_{\overline{f_1}}- [d\lfloor f_2\rfloor]_{\overline{f_2}})  \
 mod(LieR_{(\lfloor\  \rfloor, \cdot)}, w).
\end{eqnarray*}
Then
$\langle h_{f_1},h_{f_2}\rangle_w\equiv 0 \ mod(LieR_{(\lfloor\  \rfloor, \cdot)},w)$  if and only if
$$
 Red_{S_{(\cdot)}}([\lfloor f_1 \rfloor c]_{\overline{f_1}} - [d \lfloor f_2\rfloor ]_{\overline{f_2}})= \sum\alpha_i Red_{S_{(\cdot)}}( [c_i  \lfloor g_i\rfloor   d_i]_{\overline{g_i}}).
$$

Therefore, $LieR_{(\lfloor\  \rfloor, \cdot)}$ is a Gr\"{o}bner-Shirshov basis if and only if the conditions $(a)$ and $(b)$ hold.

If $LieR_{(\lfloor\  \rfloor, \cdot)}$ is a Gr\"{o}bner-Shirshov basis, then by  Composition-Diamond lemma for Lie algebras,
we have that $Irr(LieR_{(\lfloor\  \rfloor, \cdot)})=A\cup Irr(R)$ is a  linear basis of $LieE_{(\mathfrak{a}, \mathfrak{b},\lfloor\  \rfloor,\cdot)}$,  where
$
Irr(R)=\{[w]|w\in ALSW(B), w\neq c\bar{r}d, r\in R, c, d \in B^* \}.
$
Therefore,   $LieE_{(\mathfrak{a}, \mathfrak{b},\lfloor\  \rfloor, \cdot)}=\mathfrak{a}\oplus\mathfrak{b}$ as linear space with the following multiplication: for any $a,a'\in A$ and $[u], [u']\in Irr(R)$,
$$
[(a+[u ])(a'+[u'])]=[aa']_\mathfrak{a}- [u']\cdot a +  [u]\cdot a' +\sum \delta_i [u_i]+  \sum\beta_jRed_{S_{(\cdot)}}([c_j  \lfloor g_j   \rfloor    d_j]_{\overline{f_j}}),
$$
where
$
([u][u'])=\sum \delta_i [u_i]+\sum\beta_j[c_jg_jd_j]_{\overline{g_j}}
$
  in $Lie(B)$ with  each  $\delta_i, \beta_j\in k,\  [u_i]\in Irr(R),  c_j, d_j\in B^*$,    $g_j\in R$.

It is easy to see that
$\mathfrak{a}$ is an ideal of $LieE_{(\mathfrak{a}, \mathfrak{b},\lfloor\
\rfloor,\cdot)}$. Define
$$
\pi_{  \llcorner\lrcorner,\cdot} : LieE_{(\mathfrak{a}, \mathfrak{b},\lfloor\
\rfloor,\cdot)}= \mathfrak{a}\oplus \mathfrak{b}\rightarrow \mathfrak{b},\ \tilde{a}+\tilde{b}\mapsto \tilde{b},\ \ \tilde{a}\in \mathfrak{a}, \tilde{b}\in \mathfrak{b}.
$$
It follows that
$
0\rightarrow \mathfrak{a} \overset{i_{ \llcorner\lrcorner,\cdot}}{\longrightarrow}  LieE_{(\mathfrak{a}, \mathfrak{b},\lfloor\
\rfloor,\cdot)} \overset{\pi_{ \llcorner\lrcorner,\cdot}}{\longrightarrow} \mathfrak{b}\rightarrow 0
$
is a short exact   sequence, where $ i_{ \llcorner\lrcorner,\cdot} $ is the inclusion map. Therefore,
$
  LieE_{(\mathfrak{a}, \mathfrak{b},\lfloor\  \rfloor,\cdot)}=  Lie(A\cup B| LieR_{(\lfloor\  \rfloor, \cdot)} )
$
is an extension of $\mathfrak{a}$ by $\mathfrak{b}$.
\hfill $ \square$

\ \

Let $\mathfrak{e}$ be an extension of $\mathfrak{a}$ by $\mathfrak{b}$, i.e.,
$
0\rightarrow \mathfrak{a}  \overset{i}\rightarrow \mathfrak{e}\overset{\pi}\rightarrow \mathfrak{b} \rightarrow 0
$
is  a short exact sequence, where  $i$ is the inclusion map.
For any $b\in B$, choose a $\tilde{b}\in \mathfrak{e}$ such that $\pi(\tilde{b})=b $ in $\mathfrak{b}$.
Let
$
\tilde{\theta}: Lie(A\cup B)\rightarrow \mathfrak{e}
$
be the unique homomorphism from  $Lie(A\cup  B)$ to $\mathfrak{e}$ such that $\tilde{\theta}(b)=\tilde{b}$ and $\tilde{\theta}(a)=a$ for any $a\in A, b\in B$.

\begin{lemma} \label{th3.6} Let $\mathfrak{e}$ is an extension of $\mathfrak{a}$ by $\mathfrak{b}$, i.e.,    $0\rightarrow \mathfrak{a}  \overset{i}\rightarrow \mathfrak{e}\rightarrow \mathfrak{b} \overset{\pi} \rightarrow 0$ is   a short exact sequence. Define
the  $Lie(B)$-module operation $\cdot$ by
$$
g\cdot a=[\tilde{\theta}(g)a],\  g\in Lie(B), a\in A
$$
and
$$
\lfloor \ \rfloor: R\rightarrow \mathfrak{a},\ f\mapsto \tilde{\theta}(f).
$$
Then,  $\mathfrak{a}$ is   a compatible $Lie(B)$-module and    $\lfloor\  \rfloor$  is a factor set   such that  the conditions  ($a$) and ($b$) in Lemma \ref{th3.5} hold.  Moreover,
$
\mathfrak{e}\cong  Lie(A\cup B|LieR_{(\lfloor\  \rfloor, \cdot)}).
$

\end{lemma}
{\bf Proof.}  Since $\mathfrak{a}$ is an ideal of $\mathfrak{e}$, we have
$[\tilde{\theta}(g)a]\in\mathfrak{a} $  for any $g\in Lie(B)$ and $ a\in  A$.  Thus,  the   operation $\cdot$ is  well defined.  It is easy to see that   $\mathfrak{a}$ is a   compatible $Lie(B)$-bimodule.

 For any $f\in R$, since $\pi(\widetilde{\theta}(f))=f =0$ in $\mathfrak{b}$,  we have
$\widetilde{\theta}(f)\in \mathfrak{a}$. Thus,  $\lfloor \ \rfloor$ is well defined. Let  $(vaw)\in (A\cup B)^{**},  v,w \in B^*, a\in A$. Then,  it is easy to see that
$
Red_{S_{(\cdot)}}( (v a  w ) )= \tilde{\theta}{((vaw))}.
$
Therefore, we have

$(a)$  For any $a\in A,  f\in R$,
$
f\cdot a =[\tilde{\theta}(f)a]=[\lfloor f\rfloor a]_{\mathfrak{a}}.
$

$(b)$ For any intersection  composition $\langle f_1, f_2 \rangle_w$, where $ f_1, f_2\in R,  w=\overline{f_1}c=d\overline{f_2}$, $c, d \in B^* $ and $|\overline{f_1}|+|\overline{f_2}|>|w|$,   if
$$
\langle f_1,f_2\rangle_w=[f_1c]_{\overline{f_1}}-[df_2]_{\overline{f_2}}=\sum\alpha_i[c_ig_id_i]_{\overline{g_i}},
$$
where each $\alpha_i\in k$, $g_i\in R$ and  $c_i\overline{g_i}d_i<w$,
then
\begin{eqnarray*}
0&=&\tilde{\theta}([f_1c]_{\overline{f_1}}-[df_2]_{\overline{f_2}}-\sum\alpha_i[c_ig_id_i]_{\overline{g_i}} )\\
&=&   \tilde{\theta}([\lfloor f_1\rfloor c]_{\overline{f_1}}-[d\lfloor f_2\rfloor ]_{\overline{f_2}}-\sum\alpha_i[c_i\lfloor g_i\rfloor d_i]_{\overline{g_i}} )\\
&=&
Red_{S_{(\cdot)}}( [\lfloor f_1\rfloor c]_{\overline{f_1}}-[d \lfloor f_2\rfloor ]_{\overline{f_2}}) -\sum\alpha_iRed_{S_{(\cdot)}}([c_i  \lfloor g_i\rfloor d_i]_{\overline{g_i}}).
\end{eqnarray*}

  Let $\tilde{\theta}(Irr(R))=\{  \tilde{\theta}([u])|[u]\in Irr(R) \}$ and $ \widetilde{\mathfrak{b}}$   the subspace of $\mathfrak{e}$ spanned by $\tilde{\theta}(Irr(R))$.  For any $[u ],[u']\in Irr(R)$, if $
([u ] [u'])=\sum \delta_i [u_i]+\sum\beta_j[c_jg_jd_j]_{\overline{g_j}}
$ in $Lie(B)$, where
each
 $\delta_i, \beta_j\in k$, $c_j, d_j\in B^*$,   $[u_i]\in Irr(R)$, $g_j\in R$, then
\begin{eqnarray*}
 [\tilde{\theta}([u])\tilde{\theta}([u'])]=\tilde{\theta}(([u]  [u']))
 &=& \sum \delta_i \tilde{\theta}([u_i])+\sum\beta_j\tilde{\theta}([c_jg_jd_j]_{\overline{g_j}})\\
 &=& \sum \delta_i \tilde{\theta}([u_i])+ \sum\beta_j Red_{S_{(\cdot)}}( [c_j\lfloor g_j\rfloor d_j]_{\overline{g_j}}).
 \end{eqnarray*}

Therefore,
$\mathfrak{e}=\mathfrak{a}\oplus \widetilde{\mathfrak{b}}$ as linear space with the multiplication: for $a, a'\in A$ and $[u ],[u']\in Irr(R)$,
\begin{eqnarray*}
&&[(a+\tilde{\theta}([u]))( a'+\tilde{\theta}([u']))]\\&=&[aa']+ [a\tilde{\theta}([u'])]+[\tilde{\theta}([u])a']  + [\tilde{\theta}([u])\tilde{\theta}([u'])]\\
&=&[aa']_{\mathfrak{a}}- a\cdot[u']+ [u]\cdot a'+\sum \delta_i \tilde{\theta}([u_i])+  \sum\beta_jRed_{S_{(\cdot)}}([c_j \lfloor g_j \rfloor d_j]_{\overline{g_j}}).
\end{eqnarray*}

Define a linear map
$$
\varphi: LieE_{(\mathfrak{a}, \mathfrak{b},\lfloor\  \rfloor, \cdot)}=\mathfrak{a}\oplus \mathfrak{b}\rightarrow  \mathfrak{e}=\mathfrak{a}\oplus \widetilde{\mathfrak{b}}, \
 a\mapsto a,\ [u]\mapsto  \tilde{\theta}([u]),\ a\in A,\ [u]\in Irr(R).
$$
Then,  it is easy to see that  $\varphi$ is an isomorphism  by Lemma  \ref{th3.5}.\hfill$ \square$\\

\noindent {\bf Remark:} By  Lemma \ref{th3.6},  for a given   extension  $0\rightarrow \mathfrak{a}  \overset{i}\rightarrow \mathfrak{e}\overset{\pi}\rightarrow \mathfrak{b} \rightarrow 0$, it is   equivalent to some extension $
0\rightarrow \mathfrak{a} \overset{i_{ \llcorner\lrcorner,\cdot}}{\longrightarrow}  LieE_{(\mathfrak{a}, \mathfrak{b},\lfloor\
\rfloor,\cdot)} \overset{\pi_{ \llcorner\lrcorner,\cdot}}{\longrightarrow} \mathfrak{b}\rightarrow 0,
$
where  $\cdot$ is a compatible $Lie(B)$-module operation on $\mathfrak{a}$ and  $\lfloor\
\rfloor$ is  a factor set such that the conditions  ($a$) and ($b$) in Lemma \ref{th3.5} hold. If this the case, we call $(\lfloor\  \rfloor, \cdot)$  a special extension pair associated to the extension  $0\rightarrow \mathfrak{a}  \overset{i}\rightarrow \mathfrak{e}\overset{\pi}\rightarrow \mathfrak{b} \rightarrow 0$.\\

Let $\eta: B\rightarrow \mathfrak{a}$ be  a map and
$\Gamma_{\eta}:Lie(A\cup B)\rightarrow Lie(A\cup B)$
be the unique homomorphism   such that $\Gamma_{\eta}(b)=b-\eta(b)$ and $\Gamma_{\eta}(a)=a$ for any $b\in B, a\in A$.   If $LieR_{(\lfloor\  \rfloor', \cdot')}$ is  a Gr\"{o}bner-Shirshov basis in
$Lie(A\cup B) $, then for any $g\in Lie( B)  \subseteq Lie(A\cup B)$, we  have $Red_{LieR_{(\lfloor\  \rfloor', \cdot')}} ( g-\Gamma_{ \eta}(g)) \in \mathfrak{a} $. \\

By   Lemmas \ref{th3.5} and  \ref{th3.6}, we have the following theorem which gives a complete characterization of   extensions of     $\mathfrak{a}$ by   $\mathfrak{b}$, where  $\mathfrak{b}$ is presented by generators and relations.

\begin{theorem} \label{th3.7}
Let $\mathfrak{a}$   be a Lie algebra with    a linear basis  $A $ and   the   multiplication table:
$[aa']=[aa']_\mathfrak{a}, a,a'\in A $ and $\mathfrak{b}=Lie(B|R)$, where   $R$ is a minimal   Gr\"{o}bner-Shirshov basis in $Lie(B)$. Then,  the following two statements hold.

\begin{enumerate}
\item[(I)]
A Lie algebra $\mathfrak{e}$
is an extension of $\mathfrak{a}$ by $\mathfrak{b}$ if and only if
$$
\mathfrak{e}\cong  LieE_{(\mathfrak{a}, \mathfrak{b},\lfloor\  \rfloor, \cdot)}= Lie\left(
  A\cup B
 \left|
 \begin{array}{ll}
\  [aa']-[aa']_\mathfrak{a},  & a, a'\in A  \\
\  [ba]-b\cdot a,     & b\in B, a\in A \\
\   f- \lfloor f \rfloor, & f \in R
 \end{array}
\right. \right)
$$
for some   compatible $Lie(B)$-module structure in $\mathfrak{a}$ (denote the module operation by $\cdot$)  and  some factor set
$\lfloor\  \rfloor : R \rightarrow \mathfrak{a}$  such that
\begin{enumerate}
\item[(a)]  For any $a\in A,  f\in R$,   $f\cdot a=[\lfloor f\rfloor a]_{\mathfrak{a}}$   in  $\mathfrak{a}$.
\item[(b)] For any intersection  composition $\langle f_1, f_2 \rangle_w$, where $ f_1, f_2\in R,  w=\overline{f_1}c=d\overline{f_2}$, $c, d \in B^* $ and $|\overline{f_1}|+|\overline{f_2}|>|w|$,   if
$$
\langle f_1,f_2\rangle_w=[f_1c]_{\overline{f_1}}-[df_2]_{\overline{f_2}}=\sum\alpha_i[c_ig_id_i]_{\overline{g_i}},
$$
where each $\alpha_i\in k$, $c_i, d_i\in B^*$,  $g_i\in R$ and  $c_i\overline{g_i}d_i<w$,
then in $\mathfrak{a}$,
$$
 Red_{S_{(\cdot)}}( [\lfloor f_1\rfloor c]_{\overline{f_1}}-[d \lfloor f_2\rfloor ]_{\overline{f_2}}) =\sum\alpha_i Red_{S_{(\cdot)}}([c_i \lfloor g_i\rfloor   d_i]_{\overline{g_i}}).
$$
\end{enumerate}

 \item[(II)]Two
  extensions $\mathfrak{e}$ and $\mathfrak{e'} $ of $\mathfrak{a}$ by $\mathfrak{b} $ are  equivalent     if and  only if  there exists a   map  $\eta:B\rightarrow
\mathfrak{a}$  such that
\begin{eqnarray*}
 b\cdot' a-b \cdot a &=& [\eta(b)a]_{\mathfrak{a}},  \\
\lfloor f\rfloor' -\lfloor f\rfloor&= &Red_{LieR_{(\lfloor \rfloor', \cdot')}}(f- \Gamma_{ \eta}(f))   \label{c2}
\end{eqnarray*}
in $\mathfrak{a}$  for any $ a\in A,\ b\in B $ and  $f\in R$,  where $(\lfloor \ \rfloor, \cdot)\ (resp.  \ (\lfloor\ \rfloor', \cdot'))$ is a special extension pair  associated to $e$ (resp. $e'$).
\end{enumerate}
\end{theorem}

{\bf Proof.}   (I) This part follows from  Lemmas \ref{th3.5} and   \ref{th3.6}.

(II)  Similar to the proof of Theorem \ref{th2.5}  (II). \hfill $\square$

\subsection{Some special cases of extensions of Lie algebras}

Let $\mathfrak{b}$ be a Lie algebra with a well-ordered linear basis $B$ and the multiplication table:   $[bb']-[bb']_\mathfrak{b}, \ b,b'\in B$.  It is easy to see that   $R=\{[bb']-[bb']_\mathfrak{b}|  b,b'\in B, b>b'\}$ is a minimal Gr\"{o}bner-Shirshov basis in $Lie(B)$ with respect to the deg-lex order on $B^*$.

Define
$
(b,b')=\lfloor [bb']-[bb']_\mathfrak{b} \rfloor, b,b'\in B,
$
which can be extended  to be a  skew-symmetric  bilinear map
$
(,):\mathfrak{b}\times \mathfrak{b} \rightarrow \mathfrak{a}.
$

Therefore, we have the following corollary  by   Theorem \ref{th3.7} and a direct  computation.
\begin{corollary}\label{co3.10} (\cite{mi})
Let $\mathfrak{a}$ (resp. $\mathfrak{b}$)  be a Lie algebra with a linear basis $A$ (resp. $B$) and  the   multiplication table:
$[aa']=[aa']_\mathfrak{a}, a,a'\in A $ (resp. $[bb']=[bb']_\mathfrak{b}, \ b,b'\in B$).   Then, a Lie algebra $\mathfrak{e}$ is an extension of $\mathfrak{a}$ by $\mathfrak{b}$
if and only if
$$
\mathfrak{e}\cong   Lie\left(
  A\cup B
 \left|
 \begin{array}{ll}
\  [aa']-[aa']_\mathfrak{a},  & a, a'\in A,    \\
\  [ba]-b\cdot a,     & b\in B, a\in A \\
\   [bb']-[bb']_\mathfrak{b}-(b,b'),    & b,b' \in B
 \end{array}
\right. \right)
$$
for some   compatible $Lie(B)$-module structure in $\mathfrak{a}$ (denote the module operation by $\cdot$)  and  some  skew-symmetric  bilinear map
  $(,):\mathfrak{b}\times \mathfrak{b} \rightarrow \mathfrak{a}$   such that for any $ b,b', b''\in B$ and $a\in A$,
$$
[bb']\cdot a- [bb']_\mathfrak{b}\cdot a =  [(b,b')a]_{\mathfrak{a}},
$$
$$
  b\cdot (b', b'')-([bb']_\mathfrak{b}, b'')+ b'\cdot (b'', b)-([b'b'']_\mathfrak{b}, b)+b''\cdot (b, b')-([b''b]_\mathfrak{b}, b')=0.
$$

\end{corollary}

\begin{corollary} \label{co3.11} Let $\mathfrak{a}$ be an  abelian Lie algebra with  a  linear  basis $A$ and $\mathfrak{b}=Lie(B|R)$, where   $R$ is a minimal   Gr\"{o}bner-Shirshov basis in $Lie(B)$.
Then  a Lie algebra   $\mathfrak{e}$
is an extension of $\mathfrak{a}$ by $\mathfrak{b}$ if and only if
$$
\mathfrak{e} \cong LieE_{(\mathfrak{a}, \mathfrak{b},\lfloor\  \rfloor, \bullet)}=  Lie\left(
  A\cup B
 \left|
 \begin{array}{ll}
\  [aa'],  & a, a'\in A,    \\
\  [ba]-b \bullet a,     & b\in B, a\in A \\
\   f-\lfloor f\rfloor,    & f \in R
 \end{array}
\right. \right)
$$
for some  $\mathfrak{b}$-module structure in $\mathfrak{a}$ (denote the module operation by $\bullet$) and some  factor set    $\lfloor\ \rfloor:R\rightarrow \mathfrak{a}$ of $\mathfrak{b}$ in $\mathfrak{a}$  such that
\begin{enumerate}
\item[(b)] For any intersection  composition $\langle f_1, f_2 \rangle_w$, where $ f_1, f_2\in R,  w=\overline{f_1}c=d\overline{f_2}$, $c, d \in B^* $ and $|\overline{f_1}|+|\overline{f_2}|>|w|$,   if
$$
\langle f_1,f_2\rangle_w=[f_1c]_{\overline{f_1}}-[df_2]_{\overline{f_2}}=\sum\alpha_i[c_ig_id_i]_{\overline{g_i}},
$$
where each $\alpha_i\in k, \ c_i, d_i\in B^*$, $g_i\in R$ and  $c_i\overline{g_i}d_i<w$,
then
$$
Red_{S_{(\bullet)}}(  [\lfloor f_1\rfloor   c]_{\overline{f_1}}-[d   \lfloor f_2\rfloor ]_{\overline{f_2}}) =\sum\alpha_i Red_{S_{(\bullet)}}([c_i \lfloor g_i  \rfloor d_i]_{\overline{g_i}}).
$$
\end{enumerate}

\end{corollary}

{\bf Proof.}   ($\Longrightarrow$)  By    Theorem \ref{th3.7},  we have $\mathfrak{a}$ is a compatible $Lie(B)$-module and
$
\mathfrak{e}\cong LieE_{(\mathfrak{a}, \mathfrak{b},\lfloor\  \rfloor, \cdot)}
$
for some factor set
  $\lfloor\  \rfloor : R \rightarrow \mathfrak{a}$  of $\mathfrak{b}$ in $\mathfrak{a}$  such that
the  conditions $(a)$ and $(b)$   in  Theorem \ref{th3.7}  hold.

Since $\mathfrak{a}$  is  abelian,     we have
$
f\cdot a= [\lfloor f\rfloor a]_{\mathfrak{a}}=0
$
for any $f\in R$ and $a\in A$. Thus,  $\mathfrak{a}$ is a  $\mathfrak{b}$-module with the
operation
$$
(g+Id(R))\bullet a=g\cdot a, \ g\in Lie(B), a\in A.
$$

($\Longleftarrow$) Define the   operation
$$
g\cdot a=(g+Id(R))\bullet a, \ g\in Lie(B), a\in A.
$$
Since $\mathfrak{a}$  is  abelian, we have
$$
g\cdot [aa']_{\mathfrak{a}}=0=[g\cdot a, a']_{\mathfrak{a}}+[a, g\cdot a']_{\mathfrak{a}}
$$
Thus, $\mathfrak{a}$ is a  compatible $Lie(B)$-module. Moreover, for any $f\in R$ and $a\in A$,
$$
f\cdot a= (f+Id(R))\bullet a=0=[\lfloor f\rfloor a]_{\mathfrak{a}}.
$$
Therefore, the conditions $(a)$ and $(b)$ in Theorem  \ref{th3.7} hold.  It is easy to see that $\mathfrak{e}\cong LieE_{(\mathfrak{a}, \mathfrak{b},\lfloor\  \rfloor, \cdot)}=LieE_{(\mathfrak{a}, \mathfrak{b},\lfloor\  \rfloor, \bullet)} $. It follows that  $\mathfrak{e}$ is an  extension of $\mathfrak{a}$ by $\mathfrak{b}$ by Theorem \ref{th3.7}.\hfill $\square$\\

By  Corollaries  \ref{co3.10}  and  \ref{co3.11}, we have the following corollary.

\begin{corollary}\label{th3.2}(\cite{ch48})
Let $\mathfrak{a}$ be an abelian Lie algebra with a linear basis $A$ and    $\mathfrak{b}$  be a Lie algebra with a linear basis  $B$  and   the   multiplication table: $[bb']=[bb']_\mathfrak{b}, \ b,b'\in B$.  Then
 a Lie algebra $\mathfrak{e}$ is an extension of $\mathfrak{a}$ by $\mathfrak{b}$ if and only if
 $$
\mathfrak{e} \cong  Lie\left(
  A\cup B
 \left|
 \begin{array}{ll}
\  [aa'],  & a, a'\in A,    \\
\  [ba]-b\bullet a,     & b\in B, a\in A \\
\   [bb']-[bb']_\mathfrak{b}-(b,b') ,    & b, b' \in B
 \end{array}
\right. \right)
$$
for some  $\mathfrak{b}$-module structure in $\mathfrak{a}$ (denote the module operation by $\bullet$)  and some     skew-symmetric  bilinear map
 $(,): \mathfrak{b}\times \mathfrak{b} \rightarrow \mathfrak{a}$ such that  for any $b,b',b''\in B$,
 $$
  b\bullet(b', b'')-([bb']_\mathfrak{b}, b'')+ b'\bullet(b'', b)-([b'b'']_\mathfrak{b}, b)+b''\bullet(b, b')-([b''b]_\mathfrak{b}, b')=0.
 $$
\end{corollary}

\subsection{Applications}

In this subsection, we  give some applications of the Theorem \ref{th3.7} and its corollaries.
\begin{example}
Let $\mathfrak{a} $ be a Lie algebra with a linear  basis $A $ and the multiplication table:
$
[aa']= [aa']_\mathfrak{a}, a,  a'\in A,
$
and $\mathfrak{b}=Lie(B|f=0)$, a  one relation Lie algebra,
where $B$ is  a set and  $f\in Lie(B)$.  Then   a  Lie algebra   $\mathfrak{e}$ is an extension of $\mathfrak{a}$ by $\mathfrak{b}$
if and only if
$$
\mathfrak{e}\cong  Lie  \left(
 A\cup B
 \left|
 \begin{array}{ll}
\ [aa']-[aa']_\mathfrak{a},  & a, a'\in A, \\
 \ [ba]-b\cdot a, &  b\in B, a\in A\\
\ f-\lfloor f\rfloor, &
 \end{array}
\right. \right )
$$
for some compatible $Lie(B)$-module structure in  $\mathfrak{a} $  and some
$\lfloor f \rfloor \in \mathfrak{a}$   such that for any  \ $a\in A$,
$
f\cdot a=[\lfloor f\rfloor a]_\mathfrak{a}.
$

Especially, if $\mathfrak{a}$ is abelian, then $\mathfrak{e}$ is an extension of $\mathfrak{a}$ by $\mathfrak{b}$
if and only if
$$
\mathfrak{e}\cong  Lie  \left(
 A\cup B
 \left|
 \begin{array}{ll}
\ [aa'],  & a, a'\in A,  \\
 \ [ba]-b\bullet a, &  b\in B, a\in A\\
\ f-\lfloor f\rfloor, &
 \end{array}
\right. \right )
$$
for some  $\mathfrak{b}$-module structure in  $\mathfrak{a}$ and some
$\lfloor f \rfloor \in \mathfrak{a}$.
\end{example}
{\bf Proof.} Noting that there is no composition  in $R=\{f\}$,  $R=\{f\}$ is a  Gr\"{o}bner-Shirshov basis in $Lie(B)$.  By Theorem \ref{th3.7}, we can obtain the result. \hfill $\square$

\begin{example}
Let $\mathfrak{a} $ be a Lie algebra with a linear basis $A $ and the multiplication table:
$
[aa']= [aa']_\mathfrak{a}, a,  a'\in A,
$
and
$
\mathfrak{b}=\mathfrak{sl}(2, k)=Lie(B|R),
$
the three dimensional simple Lie algebra, where $B=\{x, y ,h\}$ and
$
R=\{[xy]-h, [hx]-2x, [hy]+2y\}.
$
  Then a  Lie algebra    $\mathfrak{e}$ is an extension of $\mathfrak{a}$ by $\mathfrak{b}$ if and only if
$$
\mathfrak{e}\cong Lie  \left(
 A\cup B
 \left|
 \begin{array}{ll}
\ [aa']-[aa']_\mathfrak{a}, \   a, a'\in A,  \\
 \
 [xa]-x\cdot a, [ya]-y\cdot a, [ha]-h\cdot a,  \    a\in A\\
\  [xy]-h-(x,y), \ [hx]-2x-(h,x),  [hy]+2y-(h,y)\\
\end{array}
\right. \right )
$$
for  some   compatible $Lie(B)$-module structure in $\mathfrak{a} $  and some  skew-symmetric  bilinear map
  $(,):\mathfrak{b}\times \mathfrak{b} \rightarrow \mathfrak{a}$
   such that for any $a\in A$,
\begin{eqnarray*}
&&[xy]\cdot a -h \cdot a =  [(x,y) a]_\mathfrak{a}, \\
&&[hx]\cdot a-2x\cdot a=[(h,x)a ]_\mathfrak{a},\\
&&[hy]\cdot a +2y\cdot a = [(h,y)a]_\mathfrak{a},\\
&&h\cdot (x,y) + x\cdot (y,h) +y \cdot (h, x)=0.
\end{eqnarray*}
\end{example}
{\bf Proof.}  By Corollary \ref{co3.10},   we can get the result.\hfill $\square$

\end{document}